\documentclass[reqno,11pt]{amsart}
\usepackage{amscd,amssymb,verbatim}

\numberwithin{equation}{section}
\theoremstyle{plain}

\theoremstyle{plain}
\newtheorem{Thm}[subsection]{Theorem}
\newtheorem{Cor}[subsection]{Corollary}
\newtheorem{Lem}[subsection]{Lemma}
\newtheorem{Prop}[subsection]{Proposition}
\newtheorem{Conjec}[subsection]{Conjecture}

\theoremstyle{definition}
\newtheorem{Def}[subsection]{Definition}

\theoremstyle{remark}

\newtheorem{rem}[subsection]{Remark}

\newcommand{\refs}[1]{Sect. ~\ref{S:#1}}
\newcommand{\refss}[1]{Sect.~\ref{SS:#1}}
\newcommand{\reft}[1]{Theorem ~\ref{T:#1}}
\newcommand{\refl}[1]{Lemma ~\ref{L:#1}}
\newcommand{\refp}[1]{Proposition ~\ref{P:#1}}
\newcommand{\refc}[1]{Corollary ~\ref{C:#1}}

\newcommand{\refd}[1]{Definition ~\ref{D:#1}}
\newcommand{\refr}[1]{Remark ~\ref{R:#1}}
\newcommand{\refe}[1]{\eqref{E:#1}}

\newenvironment{thm}%
          { \begin{Thm}  }%
          { \end{Thm} }

\newenvironment{lem}%
          { \begin{Lem}    }%
          { \end{Lem} }

\newenvironment{propos}%
          { \begin{Prop}  }%
          { \end{Prop} }

\newenvironment{cor}%
          { \begin{Cor} }%
          { \end{Cor} }

\newenvironment{defe}%
          { \begin{Def} }%
          { \end{Def} }

\newcommand{\loc}{\operatorname{loc}}
\newcommand{\comp}{\operatorname{comp}}
\newcommand{\ecomp}{C_{c}^{\infty}(E)}
\newcommand{\fcomp}{C_{c}^{\infty}(F)}

\newcommand{\lloc}{L_{\loc}}

\newcommand{\fl}{L^{2}(F)}

\newcommand{\symd}{\widehat{D}}
\newcommand{\symdt}{\widehat{D^*}}
\newcommand{\pin}{\phi_{k}}

\newcommand{\DT}{D^*}
\newcommand{\gtm}{g^{TM}}

\newcommand{\<}{\langle}
\renewcommand{\>}{\rangle}
\newcommand{\Label}{\label}
\newcommand{\Dom}{\operatorname{Dom}}
\newcommand{\End}{\operatorname{End}}

\newcommand{\supp}{\operatorname{supp}}

\newcommand{\IM}{\operatorname{Im}}
\newcommand{\RE}{\operatorname{Re}}
\newcommand{\Tr}{\operatorname{Tr}}

\newcommand\CC{\mathbb{C}}
\newcommand\RR{\mathbb{R}}
\newcommand\del{\delta}
\newcommand\p{\partial}
\newcommand{\sgn}{\operatorname{sgn}}

\newcommand{\mcomp}{C_{c}^{\infty}(M)}
\newcommand{\esm}{C^{\infty}(E)}

\newcommand{\lploc}{L_{\loc}^{p}}

\newcommand{\loloc}{L_{\loc}^{1}}

\newcommand{\lolocv}{(L_{\loc}^{1}(\RR^n))^m}
\newcommand{\wotcomp}{W^{1,2}_{\comp}}

\newcommand{\lapg}{\Delta_M}
\newcommand{\vp}{V_{+}}
\newcommand{\vm}{V_{-}}

\newcommand{\rncomp}{C_{c}^{\infty}(\RR^n)}

\newcommand{\pari}{\partial_i}
\newcommand{\park}{\partial_k}
\newcommand{\derxv}{\frac{\partial}{\partial x^i}}
\newcommand{\deryv}{\frac{\partial}{\partial y^i}}
\newcommand{\Realna}{\operatorname{Re}}

\newcommand{\znak}{\operatorname{sign}}
\newcommand{\dist}{\operatorname{dist}}

\newcommand{\sh}{\mathcal{H}}

\newcommand{\sd}{\mathcal{D}}
\newcommand{\sj}{\mathcal{J}^{\rho}}

\newcommand{\ur}{u^{\rho}}
\newcommand{\uep}{|u|_{\epsilon}}
\newcommand{\jr}{J^{\rho}}
\newcommand{\kjr}{j_{\rho}}
\newcommand{\kr}{K^{\rho}}
\newcommand{\kkr}{k_{\rho}}
\newcommand{\qrkato}{\frac{u^{\rho}}{|u^{\rho}|_\epsilon}}
\newcommand{\qkato}{\frac{u}{|u|_\epsilon}}
\newcommand{\hvp}{h_{\vp}}

\newcommand{\pik}{\phi_{k}}

\newcommand\tilphi{{\tilde{\phi}}}
\newcommand\tilb{{\tilde{b}}}
\newcommand\tilg{{\tilde{g}}}
\newcommand\tilmu{{\tilde{\mu}}}
\newcommand\tilE{{\tilde{E}}}
\newcommand\tilH{{\tilde{H}}}
\newcommand\hatD{{\hat{D}}}

\newcommand\tilDel{{\tilde{\Del}}}
\newcommand{\oU}{{\overline{U}}}
\newcommand{\oE}{{\overline{E}}}
\newcommand{\alp}{\alpha}
\newcommand{\bet}{\beta}
\newcommand{\obet}{{\overline{\bet}}}
\newcommand{\ob}{{\overline{b}}}

\newcommand{\oW}{{\overline{{W}}}}
\newcommand{\eps}{\epsilon}
\newcommand{\sig}{\sigma}
\newcommand{\lam}{\lambda}
\newcommand{\Del}{\Delta}
\newcommand{\Ome}{\Omega}
\newcommand{\Lam}{\Lambda}

\newcommand{\n}{\nabla}
\newcommand{\tiln}{\tilde{\n}}
\newcommand\Id{\operatorname {Id}}
\newcommand{\sign}{\operatorname{sign}}
\renewcommand{\div}{\operatorname{div}}
\newcommand{\grad}{\operatorname{grad}}
\newcommand\nek{,\ldots,}
\newcommand{\wl}{W^{1,2}_{\loc}(E)}

\def\H{{\mathcal H}}
\def\pa{\partial}
\def\R{\mathbb R}
\newcommand{\oA}{\bar{A}}
\newcommand\esssup{\rm ess\,\sup}

\begin{document}
\title{Essential self-adjointness of Schr\"odinger type operators on manifolds}
\author{Maxim Braverman, Ognjen Milatovic, Mikhail Shubin}

\dedicatory{Dedicated to M. I. Vishik on occasion of his 80th birthday}
\thanks{Research of M. Shubin was partially supported by the NSF
grant DMS 0107796}

\address{Department of Mathematics\\
           Northeastern University   \\
           Boston, MA 02115 \\
           USA
            }
\email{maxim@neu.edu; ogmilato@lynx.dac.neu.edu; shubin@neu.edu}

\begin{abstract}
We obtain several essential self-adjointness conditions for the Schr\"odinger type
operator $H_V=D^*D+V$, where $D$ is a first order elliptic differential operator acting
on the space of sections of a hermitian vector bundle $E$ over a manifold $M$ with
positive smooth measure $d\mu$, and $V$ is a Hermitian bundle endomorphism. These
conditions are expressed in terms of completeness of certain metrics on $M$ naturally
associated with $H_V$. These results generalize the theorems of E.~C.~Titchmarsh,
D.~B.~Sears, F.~S.~ Rofe-Beketov, I.~M.~Oleinik, M.~A.~Shubin and M.~Lesch. We do not
assume a priori that $M$ is endowed with a complete Riemannian metric. This allows us to
treat e.g. operators acting on bounded domains in $\RR^n$ with the Lebesgue
measure. We also allow singular potentials $V$. In particular, we obtain a new
self-adjointness condition for a Schr\"odinger operator on $\RR^n$ whose potential has the
Coulomb-type singularity and is allowed to fall off to $-\infty$ at infinity.

{}For a specific case when the principal symbol of $D^*D$ is scalar, we establish more
precise results for operators with singular potentials. The proofs are based on an
extension of the Kato inequality which modifies and improves  a result of H.~Hess,
R.~Schrader and D.~A.~Uhlenbrock.
\end{abstract}

\maketitle \setcounter{tocdepth}{1} \tableofcontents


\section{Introduction}\label{S:introd}

Spectral theory of differential operators takes its deep roots in non-relativistic
quantum mechanics, even though the first results about eigenvalues and eigenfunctions are
much older. Let us recall that according to von Neumann a state of a quantum system is a
vector $\psi$ in a complex Hilbert space $\H$ ($\psi$ is defined up to a non-vanishing
complex  factor), whereas  a physical quantity is a self-adjoint operator $A$ in $\H$. A
special  physical quantity $H$, called {\it energy} or {\it hamiltonian}, is responsible
for the evolution of any given state which is given by the Schr\"odinger equation
\begin{equation}\label{E:Schro}
        \frac{1}{i}\frac{\pa\psi}{\pa t} \ = \ H\psi,
\end{equation}
where we choose units so that the Planck constant is 1. Given an initial condition
\begin{equation}\label{E:Schro-ini}
        \psi(0) \ = \ \psi_0,
\end{equation}
we can write
\begin{equation}\label{EE:Schro-solu} \notag
        \psi(t) \ = \ \exp(itH)\,\psi_0,
\end{equation}
where the exponent can be defined by the spectral theorem, and the equation
\eqref{E:Schro} is satisfied in the strong sense if we assume that
$\psi_0\in\Dom(H)$.

Concrete hamiltonians $H$ come from a vaguely defined procedure which is applied to  a
classical system and called {\it quantization}. This procedure leads to hamiltonians $H$
which are usually second order differential operators, called Schr\"odinger operators.
However, it is a rare situation when one can actually immediately see that such an
operator $H$ is self-adjoint in a natural $L^2$ space. One can usually only see that $H$
is {\it symmetric} (or, in a slightly different terminology, {\it hermitian}). The
difference between symmetric and self-adjoint operators is often ignored by physicists,
even very good ones, but this difference is very important. Indeed, for a symmetric
operator $H$ a solution of the initial value problem
\eqref{E:Schro}-\eqref{E:Schro-ini}
may not exist, even if $\psi_0\in\Dom(H)$.  We can try to extend $H$ to make it a maximal
operator (i.e. make it defined on all $\psi\in L^2$  such that $H\psi\in L^2$ where
$H\psi$ is understood in the sense of distributions). But then it may happen that the
solution of \eqref{E:Schro}-\eqref{E:Schro-ini} is not unique.

Existence and uniqueness of the solution of
\eqref{E:Schro}-\eqref{E:Schro-ini} is
guaranteed (for $\psi_0\in\Dom(H)$) only if $H$ is self-adjoint. It is usual, however,
that one cannot immediately see what the natural domain of $H$ is. Then one starts from
smooth functions with compact support and tries to see whether $H$ can be {\it uniquely}
extended to a self-adjoint operator.  If this is the case, then $H$ is called {\it
essentially self-adjoint}. Under weak restrictions on regularity of coefficients this is
equivalent to saying that the closure of the operator $H$ from the smooth functions with
compact support coincides with the maximal operator. Then one has the only natural
self-adjoint operator associated with the given differential operator, and the evolution
of the quantum system is naturally well defined. By this reason essential
self-adjointness is often referred to as {\it quantum completeness}.

It is clear from the arguments above that establishing essential self-adjointness of a
Schr\"odinger type operator must be a starting point of any further investigation of the
corresponding quantum system, where this operator serves as a hamiltonian. Therefore it
is a fundamental problem of mathematical physics to provide precise and effective
conditions for essential self-adjointness.

The problem of essential self-adjointness for Schr\"odinger type operators was mainly
studied in $L^2(\R^n)$. There are thousands of papers devoted to such a study,
the first one (though in a different terminology) is due to H.~Weyl (1909). But it makes
perfect sense to consider the situation in a curved space or, more generally, on
manifolds, especially on Riemannian manifolds. Here one usually meets an interesting
interaction of analysis and geometry, since geometry of the manifold can play an
important role. The first paper in this direction is due to M.~Gaffney (1954), where
essential self-adjointness of the Laplacian on any complete Riemannian manifold was
proven.

In this paper we attempt an approach which unifies and extends almost all earlier results
on essential self-adjointness of Schr\"odinger type operators in vector bundles on
manifolds. The main results are of two types. An important feature of the first  type is
that we do not assume a Riemannian metric to be given a priori, but construct it from the
operator. Another important feature comes from a close watch of trajectories of the
corresponding classical hamiltonian and goes back to I.~Oleinik (1993, 1994).
Our work was mainly inspired by his results and also by recent work by M.~Lesch (2000) who
extended Oleinik's results and methods to a much more general context of
Schr\"odinger type operators in sections of vector bundles. Adding non-isotropy in
momentum variables we eventually came to first essential self-adjointness results  which
are non-isotropic both in space and momentum variables for the Lesch type operators with
singular potentials in sections of vector bundles.

The second type of results is obtained by an extension of the Kato inequality to vector
bundles which was first proclaimed by H.~Hess, R.~Schrader and D.~A.~Uhlenbrock (1980).
However they only stated their inequality for smooth sections, and this is not sufficient
for applications to essential self-adjointness if a non-smooth potential is added. We
rectify and improve their result, then apply it to transfer essential self-adjointness
results which were previously known for scalar operators only, to operators which have
scalar principal symbol.

One of the main tools which allows application of the Kato inequality is the positivity
of the Schwartz kernel of the resolvent of the Laplacian. In fact, a slightly stronger
result is needed. Such a result is obtained in Appendix~B, where we also discuss other
related questions and formulate a conjecture, which, if proven, would simplify our proof.
This conjecture also has an independent interest.

All main results of this paper are actually new not only on manifolds but in $\R^n$ as
well.

We tried to make this paper as self-contained as possible. To this end we sometimes chose
to repeat arguments of older papers, even in case where only minor modifications were
needed to obtain what was needed for our goals. We did it sometimes also for some
obviously known results if no references were easily available. However in most cases we
improved the exposition given in older papers and rendered them in a more modern
language. The most standard material is moved to appendices.

\section{Main results}\label{S:main}

\subsection{The setting}\label{SS:setting}
Let $M$ be a $C^{\infty}$-manifold without boundary, $\dim M=n$. We will always assume
that $M$ is connected.  We will also assume that we are given a positive smooth measure
$d\mu$, i.e. in any local coordinates $x^{1}, x^{2},\dots,x^{n}$ there exists a strictly
positive $C^{\infty}$-density $\rho(x)$ such that $d\mu=\rho(x)dx^{1}dx^{2}\dots dx^{n}$.
Let $E$ be any hermitian vector bundle over $M$.  We denote by $L^2(E)$ the Hilbert space
of square integrable sections of $E$ with respect to the scalar product
\begin{equation}\label{EE:inner} \notag
       (u,v) \ = \
           \int_{M}\, \< u(x),v(x)\>_{E_{x}}\, d\mu(x).
\end{equation}
Here $\<\cdot,\cdot\>_{E_{x}}$ denotes the fiberwise inner product.

Let $F$ be another hermitian vector bundle on $M$. Consider a first order
differential operator $D\colon \ecomp\to\fcomp$ (here $C^\infty_c$ denotes the space of
smooth compactly supported sections).  We assume that the principal symbol of $D$ is
injective.  In other words, $D$ is elliptic (possibly overdetermined).

Let $\DT$ be the formal adjoint of $D$, so $\DT:\fcomp\to\ecomp$ is a differential
operator such that
\[
       (Du,v) \ = \ (u,\DT v),
\]
for all $u\in\ecomp$, $v\in\fcomp.$

The main purpose of this paper is to give a sufficient condition for the essential
self-adjointness of the operator
\begin{equation}\label{E:HV}
       H_V \ = \ D^*D \ + \ V,
\end{equation}
where $V\in L^2_{\loc}(\End{E})$ is a linear self-adjoint bundle
endomorphism\footnote{The assumption that $V$ is locally square-integrable is needed in
order to have $H_Vu\in L^2_{\loc}(E)$ for $u\in C^\infty_c(E)$.}, i.e., for all $x\in M$,
the operator $V(x)\colon E_{x}\to E_{x}$ is self-adjoint.

Let us recall definitions of the minimal and maximal operators associated with the
differential expression $H_V$. The minimal operator $H_{V,\min}$ is the closure of $H_V$
from the initial domain $C^\infty_c(E)$, whereas the maximal operator can be defined as
$H_{V,\max}=(H_{V,\min})^*$. In particular, the domain $\Dom(H_{V,\max})$ of $H_{V,\max}$
coincides with the set of sections $u\in L^2(E)$, such that $H_Vu\in L^2(E)$ (here the
expression $H_Vu$ is understood in the sense of distributions). Then showing that $H_V$ is
essentially self-adjoint amounts to proving that $H_{V,\max}$ is a symmetric operator.

We make the following assumption on $V$.

\subsection*{Assumption A}\label{SS:Asa}
$V=\vp+\vm$, where
\begin{enumerate}
\item $\vp(x)\geq0, \vm(x)\leq0$ as linear operators $E_x\to E_x$ for every $x\in M$.
\item for every compact $K\subset M$ there exist positive constants $a_K<1$
      and  $C_K$ such that
\begin{equation}\label{E:V<Del}
      \Big(\, \int_K\, |V_-|^2\, |u|^2\, d\mu\, \Big)^{1/2}
      \ \le \ a_K\|\Del_M\, u\| \ + \ C_K\|u\|, \qquad
      \text{for all} \quad u\in C^\infty_c(M),
\end{equation}
where $|V_-(x)|$ denotes the norm of the linear map $V_-(x):E_x\to E_x$, $\Del_M:=d^*d$
is the scalar Laplacian of an arbitrary Riemannian metric on $M$ (cf. \refd{Bochner}),
and $\|\cdot\|$ is the norm in $L^2(E)$.
\end{enumerate}

\begin{rem}\label{R:Lp-dom}
The condition \refe{V<Del} is automatically satisfied if
$V_-\in L^p_{\loc}(\End{E})$ with $p\ge\frac{n}{2}$ for $n\ge5$,
$p>2$ for $n=4$, and $p=2$ for
$n\le 3$. (If $K$ is contained in a coordinate neighborhood this
follows from Theorem~IX.28, arguments from the proof of Theorem X.15, and
Theorems~X.20, X.21 of \cite{rs}. The
general case may be proven using a localization technique as it is explained in
\cite[\S5.2]{sh}.) Another option is to require that $V_-\in S_{n,\loc}$ where $S_{n,\loc}$ is a local Stummel
class, cf. Appendix~C.
\end{rem}

Assumption~A is closely related to the regularity of sections which belong to the maximal
domain $\Dom(H_{V,\max})$ of the operator $H_V$. More precisely, the following results
hold:
\begin{thm}\label{T:domain}
\begin{enumerate}
\item\label{domain1}
Suppose that $V\in L^p_{\loc}(\End E)$, with $p>\frac{n}{2}$, for $n\ge4$; and $p=2$, for
$n\le 3$. Then $\Dom(H_{V,\max})\subset W^{2,2}_{\loc}(E)$.
\item\label{domain2}
Suppose $V$ satisfies Assumption~A and the operator $D^*D$ has a scalar principal symbol.
Then $\Dom(H_{V,\max})\subset\wl$.
\end{enumerate}
\end{thm}
Note that due to Remark \ref{R:Lp-dom} the requirement on $V$ in \reft{domain}(\ref{domain1})
implies that the Assumption~A is satisfied for $V$.

The proof of \reft{domain}(\ref{domain1}) is quite simple and is given in
\refs{prgeneral}. The proof of \reft{domain}(\ref{domain2}) is given in \refs{domain} and
is based on the generalization of the Kato inequality which we obtain in \refs{Kato}.

Motivated by \reft{domain} we suggest the following conjecture (which is not trivial even
in case $M=\RR^n$):
\begin{Conjec}
If $V$ satisfies Assumption~A, then $\Dom(H_{V,\max})\subset\wl$.
\end{Conjec}
Note that \reft{domain}(\ref{domain1}) implies, in particular, that the conjecture holds
if $n\le 3$. Also, by \reft{domain}(\ref{domain2}), the conjecture is true if $H_V$ acts
on scalar functions (e.g., if $H_V$ is the magnetic Schr\"odinger operator).

Since we were not able to prove the conjecture in its full generality, we will work in
this paper under the additional
\subsection*{Assumption~B}
$\Dom(H_{V,\max})\subset\wl$ and $\Dom(H_{V_+,\max})\subset \wl$.

Before formulating the main result, we shall introduce some additional notation.
\subsection{The metric associated to $D$}\label{SS:Fmetric}
Let $\symd\colon T^*M\otimes E\to F$ be a morphism of vector bundles defined by
\begin{equation}\label{E:defsym}
       D(\phi u) \ = \ \symd(d\phi)u+\phi Du,
\end{equation}
where $u\in C^{\infty}(E)$, $\phi\in C^{\infty}(M)$, and $\hat{D}(d\phi)(u)$ is
identified with $\hat{D}(d\phi\otimes u)$, so $\hat{D}=-i\sigma(D)$, where $\sigma(D)$ is
the principal symbol of $D$.  Note that $\symdt(\xi)=-(\symd(\xi))^{*}$, where $\xi\in
T_{x}^{*}M$.

If $\xi\in T^*_xM$, then $\symd(\xi)$ defines a linear operator $E_x\to F_x$.  For
$\xi,\eta\in T^*_xM$, define
\begin{equation}\label{E:trace}
       \< \xi,\eta \> \ = \
       \frac1m\, \RE\, \Tr\left(\left(\symd(\xi)\right)^*\symd(\eta)\right),
       \qquad m= \dim E_x,
\end{equation}
where $\Tr$ denotes the usual trace of a linear operator.  Since $D$ is an elliptic
first-order
differential operator, $\symd(\xi)$ is linear in $\xi$, it is easily checked that
\refe{trace} defines an inner product on $T^*_xM$.  Its dual defines a Riemannian metric
on $M$ which we denote by $g^{TM}$.
\begin{rem}\label{R:equivalence}
If we define $|\xi|_0=|\symd(\xi)|$, where $|\symd(\xi)|$ denotes the usual norm of the
linear operator $\symd(\xi)\colon E_x\to F_x$, then $|\xi|_0$ induces a norm on
$T^*_{x}M$.  Its dual norm on $T_xM$ induces a Finsler metric on $M$. Elementary linear
algebra shows that for every $\xi\in T^*_xM$,
\begin{equation}\label{E:equalst}
      m^{-1/2} |\xi|_0 \ \leq \ |\xi| \ \leq \ |\xi|_0,
\end{equation}
where $|\cdot|$ is the norm defined by \refe{trace}.
\end{rem}
We say that a curve  $\gamma\colon [a,\infty)\to M$ {\em goes to infinity} if for any
compact set $K\subset M$ there exists $t_K>0$ such that $\gamma(t)\not\in K$, for all
$t\ge t_K$. The metric $g^{TM}$ is called {\em complete} if $\int_{\gamma}ds= \infty$ for
every curve $\gamma$ going to infinity. Here $ds$ is the element of arclength
corresponding to the metric $g^{TM}$.   The completeness of $g^{TM}$ is equivalent to
geodesic completeness of $M$.  In particular, completeness conditions on metrics in
\refr{equivalence} corresponding to $|\cdot|_0$ and $|\cdot|$ are equivalent.

We say that $\int^{\infty}\frac{ds}{\sqrt{q}}=\infty$ for a function $q\colon M\to\RR$
if $\int_{\gamma}\frac{ds}{\sqrt{q}}=\infty$, for every curve $\gamma$ on $M$ going to
infinity.

The main result of this paper is the following
\begin{thm}\label{T:main}
Suppose $V$ satisfies Assumptions~A and B above. Assume that $q\colon M\to\RR$, and the
following conditions are satisfied
\begin{enumerate}
\item
$q\geq1$ and $q^{-1/2}$ is globally Lipschitz, i.e. there exists a constant $L>0$ such
that, for every $x_{1}, x_{2}\in M$,
\[
       |q^{-1/2}(x_{1})-q^{-1/2}(x_{2})|\leq Ld(x_{1},x_{2}),
\]
where $d$ is the distance induced by the metric $\gtm$.
\item There exists $\delta\in[0,1)$ such that $$\delta D^*D+V\geq-q$$
(in the sense of quadratic forms on $\ecomp$)

\item $\int^{\infty}\frac{ds}{\sqrt{q}}=\infty$, where $ds$ is the
arclength element of $\gtm$.
\end{enumerate}
Then $H_V$ is essentially self-adjoint on $\ecomp$.
\end{thm}
The proof is given in \refss{prmain}.
\begin{rem}\label{R:rmk}
Since we assume $q\ge1$, condition (iii) of the theorem implies that $\gtm$ is complete.
Note also that condition (iii) itself is equivalent to the statement that the new metric
$g:=q^{-1}g^{TM}$ is complete.
\end{rem}
\smallskip
 From \reft{main}, we immediately obtain the following
\begin{cor}\label{C:sled1}
Suppose that $V\geq-q$ where $q:M\to\RR$ satisfies assumptions (i) and (iii) of the above
theorem. Then $H_V$ is essentially self-adjoint.
\end{cor}
\begin{cor}\label{C:noV}
If the metric $g^{TM}$ is complete, then the operator $D^*D$ is essentially self-adjoint.
\end{cor}
\begin{rem}
Note that both the space $L^2(E)$ and the operator $H_V$ depend on the measure $d\mu$ on
$M$. However, this measure did not appear in the formulation of \refc{sled1}. Hence, if
the conditions of this corollary are satisfied, the operator $D^*D+V$ will be essentially
self-adjoint for any choice of the measure.  In particular, $M$ may have a finite volume
with respect to this measure.
\end{rem}
\begin{rem}\label{R:Laplacian}
\reft{main} covers some interesting operators, appearing in differential geometry. In
particular, if $D=d: C^\infty(M)\to \Ome^1(M)$ is the de Rham differential,
then $D^*D=\Del_M$ is the scalar Laplacian on $M$, cf. \refd{Bochner}.
More generally, we  can consider the de Rham differential
$d:\Ome^j(M)\to \Ome^{j+1}(M)$ on forms of arbitrary degree.
Let $d^*$ be the formal adjoint of $d$ and set
$D=d+d^*$. Then
$D^*D$ is the Laplace-Beltrami operator on differential forms, cf. \cite[\S3.6]{bgv}.  Also if
$D=\n:C^\infty(E)\to \Ome^1(M,E)$ is a covariant derivative corresponding to a metric
connection on Hermitian vector bundle $E$ (\refss{cover}), then $D^*D$ is the Bochner
Laplacian on $E$, cf. \refd{Bochner}.
Note, however, that the operators above are generally different from the classical ones since
they are defined by using an arbitrary positive smooth measure $d\mu$ which does not
necessarily coincide with the standard Riemannian  measure.
\end{rem}
\smallskip

The restriction $\delta<1$ is essential in the condition (ii) of \reft{main}. In
\refss{nonsa}, we present an example of an operator, which is not essentially
self-adjoint, but satisfies the conditions of \reft{main} with any $\delta>1$.
Unfortunately, we do not know whether \reft{main} remains true if we allow $\del=1$ in
condition (ii). However, in \refs{mpf} we prove the following
\begin{thm}\label{T:bound}
Suppose that $V$ satisfies Assumptions~A and B. Assume that the metric $g^{TM}$ is
complete and that the operator $H_V=D^*D+V$ is semi-bounded below on $C^\infty_c(E)$.
Then $H_V$ is essentially self-adjoint.
\end{thm}
{}For the case when $H_V$ is a magnetic Schr\"odinger operator acting on scalar
functions, this result was formulated in~\cite{sh}. Unfortunately, in case of singular
potentials not all details of the proof were presented. These details, however, can be
found in the present paper.

\reft{main} immediately implies the following

\begin{cor}\label{C:btri}
Assume that $V\in L^2_{\loc}(\End E)$ satisfies Assumption~A. Suppose that
$V=V_{1}+V_{2}$, where $V_{1}, V_2\in L^2_{\loc}(\End E)$ are such that the operator
$\delta D^*D+V_1$ is semi-bounded below on $C_c^\infty(E)$ for some $\delta<1$, and
$V_{2}\geq -q$, where $q$ satisfies conditions (i) and (iii) of \reft{main}. Then the
operator $H_V=D^*D+V$ is essentially self-adjoint.
\end{cor}

\begin{rem}\label{R:btri}
{}If we assume that $M=\RR^{n}$ with the standard metric and measure,
$H_V=-\Delta+V$ where $\Delta$ is the
standard Laplacian and $V$ is a scalar (real-valued) function,
the following more explicit condition on
$V_1$  implies the semi-boundedness condition given
in the formulation of Corollary~\ref{C:btri}:
$V_{1}\in L^{p}(\RR^{n})$, where
$p\ge\frac{n}{2}$ for $n\ge 5$, $p>2$ for $n=4$, and $V_{1}\in L^{2}(\RR^{n})$
for $n\le 3$. (See references in Remark \ref{R:Lp-dom}.)

It is also sufficient to require that $V_1\in S_n$ (the Stummel class) or  $V_1\in K_n$
(the Kato class) -- see  Appendix C for definitions and properties of these classes.
\end{rem}

\subsection{Most recent history}\label{SS:hist}
Here we provide references for the most recent papers which
were our source of inspiration. For more history see Appendix D.

\reft{main} generalizes recent work of I. Oleinik~\cite{ol} and M. Shubin~\cite{sh1}. The
latter author considered the scalar magnetic Schr\"odinger operator $H_V=-\Delta_{A}+V$
on a complete Riemannian manifold, where $\Delta_{A}=d_{A}^{*}d_{A}$, $A$ is a real
sufficiently regular 1-form, and $d_{A}u=du+iuA$ is a deformed differential, while $V\in
L_{\loc}^{\infty}(M)$. Keeping the $L_{\loc}^{\infty}$ assumption, M. Braverman~\cite{br}
generalized the work of I. Oleinik to Schr\"odinger operators on differential forms.
Later, M.~Lesch~\cite{lm} noticed that one does not need to restrict oneself to
Laplacian-type operators with isotropic symbols.  He considered the generalized
Schr\"odinger operator $H_V=\DT D+V$ on a complete Riemannian manifold, where $D$ is the
same as in our paper.  He established the following sufficient condition for the
essential self-adjointness of $H_V$.

\begin{thm}[\textbf{M.~Lesch}]\label{T:Lesch}
Assume that $M$ is a complete Riemannian manifold. Assume that $V\in
L^{\infty}_{\loc}(\End E)$, and $V\geq-q$, where $q\geq1$ is a locally Lipschitz
function. Set \/
    \[
       c(x) \ := \
         \max\big\{\, 1,\,
           \sup\{|\symd(x,\xi)|\colon\, |\xi|_{T_{x}^{*}M}\leq1\}
           \, \big\}.
    \]
If
\begin{enumerate}
\item  there exists a constant $C>0$ such that $c(x)|dq^{-1/2}(x)|\leq
C$ for all $x\in M$ and

\item $\int^{\infty}\frac{1}{c\sqrt{q}}ds=\infty$,
\end{enumerate}
then $H_V$ is essentially self-adjoint.
\end{thm}

Note that \reft{main} does not a priori assume that $M$ is a Riemannian manifold. More
importantly, while Lesch uses the function $c(x)$ which is constructed in isotropic
fashion (by taking supremum over $|\xi|=1$, $\xi\in T^{*}M$), we make a more effective
use of the symbol $\symd$ by considering the metric $\gtm$. In \refss{ex1} we construct
an explicit operator for which \reft{main} guarantees the essential self-adjointness,
while Lesch's theorem is inconclusive.

\begin{rem}
M.~Lesch \cite{lm} considered a more general case, where the operator $D$ is not elliptic
but satisfies the following condition: \/ if $u\in\Dom(H_{V,\max})$ then
$Du\in L^2_{\loc}(F)$. (It follows from \reft{main} that this condition is automatically
satisfied if $D$ is elliptic and $V$ is sufficiently regular, e.g. $V\in L^p_{\loc}(\End E)$
where $p$ is as in \reft{domain}(i)). \reft{main} can
also be extended to non-elliptic case.
\end{rem}

\subsection{Acknowledgment}
The authors are grateful to E.~B.~Davies, A.~Grigoryan, E.~Hebey
and I.~Verbitsky for useful discussions.

\section{Some examples}

\subsection{Example: an operator with non-isotropic symbol}\label{SS:ex1}
We now consider an example of an operator with a non-isotropic symbol (a function on
$T^{*}M$ depending not only on the norm of a covector but also on its direction).  We
shall use \refc{noV} to prove that this operator is essentially self-adjoint.  This
result cannot be obtained by ``isotropic estimates" of I.~Oleinik \cite{ol}, M.~Shubin
\cite{sh1}, and M.~Lesch \cite{lm}.

Let $M=\RR^{2}$ with the standard metric and measure, and $V=0$.  Consider the operator
\[
      D \ = \ \left(\begin{array}{c}a(x,y)\frac{\partial}{\partial x}\\
                          b(x,y)\frac{\partial}{\partial y}\\
                      \end{array}\right),
\]
where $a,b$ are smooth real-valued nowhere vanishing functions in $\RR^2$.
It is elliptic (overdetermined).
We are interested in the operator
\[
      H \ := \ D^*D
      \ = \
      -\frac{\partial}{\partial x}
        \left(a^{2}\frac{\partial}{\partial x}\right)-
           \frac{\partial}{\partial y}
            \left(b^{2}\frac{\partial}{\partial y}\right).
\]

The matrix
of the inner product on $T^*M$ defined by $D$ via \refe{trace} is ${\rm diag}(a^2/2,b^2/2)$.
The matrix  of the corresponding Riemannian metric $g^{TM}$ on $M$ is
${\rm diag}(2a^{-2},2b^{-2})$, i.e. the metric itself is $ds^2=2a^{-2}dx^2+2b^{-2}dy^2$.
By \refc{noV}, to prove
that the operator $H$ is essentially self-adjoint, it is enough to show that the metric
$g^{TM}$ is complete, i.e., that
\begin{equation}\label{E:polnost}
      \int^{\infty}ds \ = \ \infty,
      \textrm{ where $ds$ is the arc-length element associated to $g^{TM}$}.
\end{equation}

Let $\gamma(t)=(x(t),y(t))$, where $t\in[0,\infty)$ is a curve in $\RR^{2}$ which goes to
infinity, cf. \refss{Fmetric}. Then the completeness condition \refe{polnost} can be
written as
$\int_{0}^{\infty}\frac{|x'(t)|}{a}dt+\int_{0}^{\infty}\frac{|y'(t)|}{b}dt=\infty$. The
purpose of this example is to show that the later condition can be satisfied even when
the integral in the left hand side of condition (ii) of \reft{Lesch} converges. Roughly
speaking, we would like to construct an example, where at least one of the integrals of
$1/a$ and $1/b$ is large, while the integral of $1/\sqrt{a^2+b^2}$ is small. This can be
achieved, for instance, as follows. Define
\begin{eqnarray}
     a(x,y) \ &=& \ (1-\cos(2\pi e^{x}))x^{2}+1;\notag\\
     b(x,y) \ &=& \ (1-\sin(2\pi e^{y}))y^{2}+1.\notag
\end{eqnarray}

Let us show that with this choice of $a,b$ the operator $H$ is essentially self-adjoint.
Let $\gamma(t)$ be as above. Then there exists a sequence of numbers $t_n\in[0,\infty)$
such that either $x(t_n)\to\infty$ or $y(t_n)\to\infty$ as $n\to\infty$.  Suppose that
$x(t_n)\to\infty$ (the other case is treated similarly). Then
\[
      \int_{0}^{t_n}\frac{|x'(t)|}{a}dt
      \ \geq \
      \int_{0}^{x(t_n)}\frac{dx}{a}.
\]
Hence, letting $n\to\infty$, we obtain
\[
      \int_{0}^{\infty}\frac{|x'(t)|}{a}dt
      \ \geq \
      \int_{0}^{\infty}\frac{dx}{a}.
\]
Thus, for \refe{polnost} to be satisfied, it is sufficient to have
$\int_{0}^{\infty}\frac{dx}{a}=\infty$.

Denote by $r_{k}$ the solutions of $\cos(2\pi e^{x})=1$, i.e. $r_{k}=\ln{k}$, $k=1,2,\dots$.
Then there is an interval $I_{k}$ around $r_{k}$ such that $1-\cos(2\pi e^{x})\leq
x^{-2}$, i.e. such that $|\sin(\pi e^x)|\leq(\sqrt{2}x)^{-1}$. Note that
\begin{equation}\label{E:1/a>}
      \frac1{a(x,y)} \ \ge \ \frac12, \qquad\text{for every}\quad
                      x\in I_k.
\end{equation}

To estimate the length of interval $I_k$ from above, notice that for every $h>0$,
\[
      |\sin(\pi e^{r_{k}+h})|
      \ = \
      |\sin(\pi e^{r_{k}+h}
           -\pi e^{r_{k}})|\leq\pi e^{r_{k}}(e^{h}-1).
\]
Thus $I_k$ includes, in particular, all $r_k+h$ with $0<h<1/2$, such that
$\pi e^{r_{k}}(e^{h}-1)<\frac{1}{\sqrt{2}(r_{k}+h)}$.  Clearly $h<r_{k}$ for $k\ge 2$.
Therefore, for $k\ge 2$, the interval $I_k$ contains all numbers $r_k+h$, such that
$0<h<1/2$ and $\pi e^{r_{k}}(e^{h}-1)<\frac{1}{2\sqrt{2}r_{k}}$. Note also
that $h<e^h-1<2h$ if $0<h<1/2$. It follows that
$I_k$ contains all $r_k+h$ with $0<h\le \frac{1}{4\sqrt{2}\pi(\ln{k})k}$.  Hence,
the length of $I_k$ is greater than $\frac{1}{4\sqrt{2}\pi(\ln{k})k}$. In view of
\refe{1/a>}, this implies that for any $x_0>0$, there exists $k_0$, such that
\[
    \int_{x_0}^{\infty}\frac{dx}{(1-\cos(2\pi e^{x}))x^{2}+1}
    \ \geq \
    \frac{1}{2}\sum_{k=k_0}^{\infty}\frac{1}{4\sqrt{2}\pi(\ln{k})k}
    \ = \ \infty.
\]
Thus by \refc{noV}, $H=D^*D$ is essentially self-adjoint.

Note that Lesch's theorem is inconclusive here.  Indeed, the Lesch's function $c$ (cf.
\reft{Lesch}) is given by
\[
      c \ =  \   \max\big\{\, 1,\,
          \sup\{|a\xi_1+b\xi_2|:\, \xi_1^2+\xi_2^2\leq1\}
          \, \big\}
      \ = \ \max\{1,\sqrt{a^2+b^2}\}.
\]
Let us consider the curve $\gamma(t)=(t,t)$. Then at a point $\gamma(t)$ we have
$a^2+b^{2}\geq(3-2\sqrt{2})t^{4}$, and condition (ii) of Lesch's theorem~\ref{T:Lesch} is
clearly violated.

\subsection{Example: an operator in a bounded domain}\label{SS:ex2}
Consider a square $S=\{(x,y)\in\RR^{2}:\, -1< x< 1, -1< y< 1\}$. Let
$D=\left(\begin{array}{c} (1-x^2)\frac{\partial}{\partial
x}\\(1-y^2)\frac{\partial}{\partial y}\end{array}\right)$. From \refc{noV} we
immediately see that the operator
\[
      H \ := \ D^*D \ = \ -\frac{\partial}{\partial x}
        \left((1-x^{2})^{2}\frac{\partial}{\partial x}\right)
          - \frac{\partial}{\partial y}
             \left((1-y^2)^{2}\frac{\partial}{\partial y}\right)
\]
is essentially self-adjoint.

\subsection{Example: Coulomb-type potential}\label{SS:ex3}
Consider $M=\RR^{3N}$, and $\mathbf{x}_{1},\dots,\mathbf{x}_{N}$ in $\RR^{3}$ orthogonal
coordinates for $\RR^{3N}$. Let
\begin{equation}\label{E:Coulomb}
    H \ = \ -\sum_{i=1}^{N}\Delta_{i} -\sum_{i=1}^{N}\frac{1}{|\mathbf{x}_{i}|}
         + \sum_{i<j}^{N}\frac{1}{|\mathbf{x}_{i}-\mathbf{x}_{j}|}  +V(x),
\end{equation}
where $\Delta_{i}$ is the Laplacian corresponding to $\mathbf{x}_{i}$, and $V(x)$ is a
locally bounded potential satisfying the following condition: \, there exists a function
$q(x)\ge 1$ such that the function $q^{-1/2}(x)$ is globally Lipschitz and
\begin{equation}\label{E:V>-q}
    V(x) \ge -q(x), \qquad\quad \int^\infty\, {ds/\sqrt{q}} \ = \ \infty.
\end{equation}
Then $H$ is essentially self-adjoint by \refc{btri} and \refr{btri}. Indeed, the terms
$|\mathbf{x}_{i}|^{-1}$ and $|\mathbf{x}_{i}-\mathbf{x}_{j}|^{-1}$ are dominated by the
Laplacian in $\R^{3N}$ both in operator sense and in the sense of quadratic forms. This
can be proved either by ``separation of variables" as in the classical Kato's paper
\cite{Kato51}, or by use of Stummel and Kato classes -- see Example~\ref{SS:Coulomb-dom}
in Appendix C for the Stummel classes (the corresponding argument can be repeated
verbatim for the Kato classes). Note that neither of these terms belongs to
$L^{3N/2}(\RR^{3N})$, except in case $N=1$.

To give a more specific example we can set $V(x)= -|x|^2$ in \refe{Coulomb} and conclude
that the operator
\[
     -\sum_{i=1}^{N}\Delta_{i}
     -\sum_{i=1}^{N}\frac{1}{|\mathbf{x}_{i}|}
         + \sum_{i<j}^{N}\frac{1}{|\mathbf{x}_{i}-\mathbf{x}_{j}|}
           -\sum_{i=1}^{N}|\mathbf{x}_{i}|^{2}.
\]
is essentially self-adjoint (This follows also from a  result of H.~Kalf \cite{Kalf73}).
More generally, the operator \refe{Coulomb} is essentially self-adjoint if there exists a
locally bounded function $r: [0,\infty)\to [1,\infty)$, such that $r^{-1/2}$ is globally
Lipschitz and
\begin{equation}\label{E:V>r}
     V(x) \ \ge \ -r(|x|), \qquad\quad
     \int_0^\infty\, {}dt/\sqrt{r(t)} \ = \ \infty.
\end{equation}
As an example, we can set $V(x)=-|x|^2\log{(1+|x|)}$.

Note that we can also consider potentials which satisfy \refe{V>-q} but not \refe{V>r}.
Constructions of rich families of such potentials can be found in \cite{RB} and
\cite{ol}. Note that self-adjointness of \refe{Coulomb} in this case can not be
established by the method of \cite{Kalf73}.

\subsection{Example: a not essentially self-adjoint operator}\label{SS:nonsa}
In this subsection we present an example of an operator, which is not essentially
self-adjoint, but satisfies the conditions of \reft{main} with every $\delta>1$.

Consider an operator $D\colon C_c^{\infty}(\RR)\to C_c^{\infty}(\RR)$, where
$D=\frac{d}{dx}+|x|^{\alpha}$, with $\alpha>3$. Then $D^*=-\frac{d}{dx}+|x|^{\alpha}$,
and hence $$D^*D=-\frac{d^{2}}{dx^{2}}-\alpha |x|^{\alpha-1} \sgn{x}+|x|^{2\alpha}.$$ Let
$V=-|x|^{2\alpha}$, and let $H_V=D^*D+V$.  Then
\[
      H_V \ = \ -\frac{d^{2}}{dx^2} \ - \ \alpha |x|^{\alpha-1}\sgn{x}.
\]
This operator {\em is not essentially self-adjoint} since $\alpha-1>2$ (cf. Example~1.1
in section~3.1 of \cite{bs}). However, for any $\delta >1$,
\[
      \delta D^*D+V
      \ = \
      -\delta \frac{d^{2}}{dx^{2}}
        \ - \ \delta\alpha |x|^{\alpha-1}\sgn{x} \ + \ (\delta-1)|x|^{2\alpha}.
\]
Since $2\alpha>\alpha-1$ and $\delta-1>0$, there exists $C_\delta>0$ such that
$-\delta\alpha{}|x|^{\alpha-1}\sgn{x}+(\delta-1)|x|^{2\alpha} > -C_\delta$ for all
$x\in\RR$. Hence, $\delta{}D^*D+V> -C_\delta$

We conclude that $\delta{}D^*D+V$ is semi-bounded below for every $\delta>1$. Therefore,
\reft{main} does not hold with $\del>1$.

\section{Proof of the first part of \reft{domain}}\label{S:prgeneral}

In this section assumptions on the potential are as in hypotheses of
\reft{domain}(\ref{domain1}).

\begin{lem}\label{L:VuL2}
{}For any $u\in \Dom(H_{V,\max})$, we have $Vu\in L^2_{\loc}(E)$.
\end{lem}
\begin{proof}
Let us assume that $u\in \Dom(H_{V,\max})$, i.e. $u\in L^2(E)$ and
\begin{equation}\label{E:eq-max}
     D^*Du+Vu \ = \ f\in L^2(E).
\end{equation}
Define $t_1:=2$. Since $V\in L^p_{\loc}(\End{E})$, the H\"older inequality implies $Vu\in
L^{s_1}_{{\loc}}(E)$, where
\[
     \frac{1}{s_1} \ = \
        \frac{1}{2}+\frac{1}{p} \ = \ \frac{1}{t_1}+\frac{1}{p}\;.
\]
Clearly, $1\le s_1<2$. Moreover, if $n>3$, then $p>2$ and  $s_1>1$.

Now we can improve $s_1$ as follows. Let us observe that
\eqref{E:eq-max} implies $D^*Du=f-Vu=h\in L^{s_1}_{\loc}(E)$. If
$s_1>1$, then by the standard elliptic regularity results (cf. \cite{Triebel}, Sect. 6.5)
we obtain $u\in W^{2,s_1}_{\loc}(E)$. Therefore by the Sobolev embedding theorem (cf.
\cite[Th.~5.4]{Adams} or \cite[Th.~4.5.8]{Horm1}) we conclude that
$u\in L^{t_2}_{\loc}(E)$, where
\[
     \frac{1}{t_2} \ = \ \frac{1}{s_1}-\frac{2}{n}
     \ = \ \frac{1}{t_1}+\frac{1}{p}-\frac{2}{n},
\]
and then by the H\"older inequality $Vu\in L^{s_2}_{\loc}$, where
\[
     \frac{1}{s_2} \ = \ \frac{1}{t_2}+\frac{1}{p} \ = \
     \frac{1}{s_1}-\frac{2}{n}+\frac{1}{p} \ = \
     \frac{1}{s_1}-\frac{2}{n}\left(1-\frac{n}{2p}\right).
\]
In this way we can continue to obtain a series of inclusions $u\in L^{t_k}_{\loc}(E)$,
$Vu\in L^{s_k}_{\loc}(E)$, $k=1,2,\dots,$ with
\[
     \frac{1}{t_{k+1}} \ = \ \frac{1}{s_k}-\frac{2}{n}, \qquad
     \frac{1}{s_{k+1}} \ = \
\frac{1}{s_k}-\frac{2}{n}\left(1-\frac{n}{2p}\right),
\]
until at some point we obtain that $s_{k}>2$ which makes the next step impossible due to
the term $f\in L^2(E)$ in \eqref{E:eq-max}. Then we can conclude that $Vu\in
L^2_{\loc}(E)$ as required.

Now assume that $n\le 3$, so $p=2$ and at the initial step we only have $Vu\in
L^1_{\loc}(E)$. By another Sobolev embedding theorem (cf. \cite[Case C in Theorem
5.4]{Adams}) we have $W^{1,q}_{\loc}(E)\subset C(E)$ provided $q>n$ (here $C(E)$ denotes
the set of continuous sections of $E$). Using standard duality arguments (cf. \cite[Sect.
3.6--3.13]{Adams}) we obtain $L^1_{\loc}(E)\subset W^{-1,q'}_{\loc}(E)$ if
$1<q'<n/(n-1)$.
Therefore, $D^*Du\in W^{-1,\tilde{q}}_{{\loc}}(E)$, where $\tilde{q}= \min\{q',2\}$. By
standard elliptic regularity results (cf. \cite[Ch. 6.5]{Triebel}) we obtain $u\in
W^{1,\tilde{q}}_{{\loc}}(E)$. In particular, $u\in L^{\tilde{q}}_{\loc}(E)$ for such
$\tilde{q}$. Since $\tilde{q}>1$, we can argue as in the first part of the proof to come
to the same conclusion that $Vu\in L^2_{\loc}(E)$.
\footnote{Note that to apply the elliptic regularity to the equation
$D^*Du=f\in W^{-1,t}_{\loc}(E)$ we must have $t>1$. If $t=1$, then the last equation does
not imply $u\in W^{1,t}_{{\loc}}(E)$, cf. \cite[Th.~7.9.8]{Horm1}. This is also the
reason we need to start the proof of the lemma with showing that $Vu\in L^t_{\loc}(E)$
with $t>1$.}
\end{proof}

\subsection{Proof of \reft{domain}(i)}

If $u\in D(H_{V,\max})$, then $D^*Du=H_Vu-Vu\in\lloc^{2}(E)$.  By elliptic regularity,
$u\in W^{2,2}_{\loc}(E)$.
\hfill$\square$

\section{Kato inequality for Bochner Laplacian}\label{S:Kato}

In this section we prove a generalization of the Kato inequality to the Bochner Laplacian
on a manifold $M$ endowed with a Riemannian metric $g=g^{TM}$ and a measure
$d\mu=\rho{}dx$. Here, $dx$ denotes the Riemannian volume form on $M$ and
$\rho:M\to(0,\infty)$ is a smooth function.

As a first step, we establish a ``smooth version" of the Kato inequality, cf.
\refp{epsmoothkato}. In the case when $d\mu$ is the Riemannian volume form on $M$, a
similar inequality was proven in H.~Hess, R.~Schrader and D.~A.~Uhlenbrock~\cite{hsu2}.
Then, in \refss{prkato} we prove ``$\loloc$" version of the Kato inequality.

\subsection{A pairing on the space of bundle-valued forms}\label{SS:bvforms}
Let $E$ be a Hermitian vector bundle over $M$ and let  $\oE$ denote the complex conjugate
of the bundle $E$. It is identical with $E$ but with multiplication by $\lambda\in\CC$
defined as multiplication by $\bar\lambda$. The identity map $E\to\oE$
is called complex conjugation and defines an anti-linear
isomorphism $E\overset{~}\to \oE$. The image of $v\in E_x$ under this isomorphism is
denoted $\bar v$. The Hermitian structure on $E$ defines a complex linear map
\begin{equation}\label{E:<>}
      \<~\>:E\otimes\oE \ \to \ \CC, \qquad a\otimes\ob \ \mapsto
        \<a\otimes\ob\> \ := \ \<a,b\>.
\end{equation}
Let $\Lam^i=T^*M\wedge\dots\wedge{}T^*M$ be the $i$-th exterior power of the cotangent
bundle to $M$. The space of smooth sections of the tensor product $\Lam^i\otimes{E}$ is
called the space of $E$-valued differential $i$-forms on $M$ and is denoted by
$\Ome^i(M,E)$. The map \refe{<>} extends naturally to the maps
\[
      \<~\>:\Lam^k\otimes E\otimes\oE\to \Lam^k,  \qquad
      \<~\>:\Ome^k(M,E\otimes\oE)\to \Ome^k(M,\CC).
\]
If $\alp\in \Ome^i(M,E), \bet\in \Ome^j(M,E)$, then $\obet\in \Ome^j(M,\oE)$,
$\alp\wedge\obet\in \Ome^{i+j}(M,E\otimes\oE)$, and $\<\alp\wedge\obet\>\in
\Ome^{i+j}(M,\CC)$. If one of the forms $\alp, \bet$ belongs to $\Ome^0(M,E)\equiv
C^\infty(E)$, then we will omit ``$\wedge$" from the notation and write simply
$\<\alp\obet\>$.

Let $*:\Lam^{i}\to \Lam^{n-i}$ be the Hodge-star operator, cf. \cite{Warner}. This
operator extends naturally to the spaces $\Lam^i\otimes{E}$ and $\Ome^i(M,E)$.

The formula
\[
      \<\alp,\bet\>_{{}_{\Lam^i\otimes{E}}} \ := \ *^{-1}\<\alp\wedge*\obet\>
      \ \in \ \Lam^0\cong\CC,
      \qquad \alp,\bet\in \Lam^i\otimes E, \ i=0\nek{n}
\]
defines a non-degenerate Hermitian scalar product on $\Lam^i\otimes{}E$, and we denote
\[
      |\alp| \ := \ \<\alp,\alp\>^{1/2}_{{}_{\Lam^i\otimes{E}}}.
\]
Note that if $\alp, \bet\in \Lam^0\otimes{}E\cong E$, then
$\<\alp,\bet\>_{{}_{\Lam^i\otimes{E}}}$ coincides with the original Hermitian scalar
product $\<\alp,\bet\>$ on $E$.

Similarly, for $\alp,\bet\in \Ome^i(M,E)$, one defines
$\<\alp,\bet\>_{{}_{\Lam^i\otimes{E}}}\in C^\infty(M)$ and $|\alp|\in C(M)$.

\subsection{Covariant derivative and its dual}\label{SS:cover}
Let $\n$ be a connection on $E$. It defines a linear map $\n:\Ome^*(M,E)\to
\Ome^{*+1}(M,E)$ such that
\[
    \n\, (\omega\wedge\bet) \ = \ d\omega\wedge \bet + (-1)^i\omega\wedge \n\bet,
    \qquad \omega\in \Ome^i(M), \ \bet\in \Ome^j(M,E).
\]
In this paper we will always assume that $\n$ is a Hermitian connection, i.e., that the
following equality is satisfied
\begin{equation}\label{E:Herm}
      d\<\alp\wedge\obet\> \ = \ \<\n\alp\wedge\obet\> \ + \
      (-1)^i\<\alp\wedge\overline{\n\bet}\>,
      \qquad \alp\in \Ome^i(M,E), \, \bet\in \Ome^j(M,E),
\end{equation}
where $d:\Ome^*(M,\CC)\to \Ome^{*+1}(M,\CC)$ is the de Rham differential.

Let $\Ome^i_{\comp}(M,E), \, i=0,1\nek n$ denote the space of compactly supported
differential forms. Consider the scalar product on the space $\Ome^i_{\comp}(M,E)$
defined by the formula
\begin{equation}\label{EE:product} \notag
      (\alp,\bet) \ = \ \int_M\, \<\alp,\bet\>_{{}_{\Lam^i\otimes{E}}}\, d\mu
      \ = \ \int_M\, \<\alp\wedge*\obet\>\, \rho, \qquad \alp,\bet\in
        \Ome^i_{\comp}(M,E).
\end{equation}
Let $\n^*:\Ome^*(M,E)\to \Ome^{*-1}(M,E)$ denote the formal adjoint of $\n$, i.e., we
have $(\n\alp,\bet)= (\alp,\n^*\bet)$ for all $\alp,\bet\in \Ome^*_{\comp}(M,E)$.
\begin{lem}\label{L:nstar}
Suppose $\bet\in \Ome^j(M,E)$. Then
\begin{equation}\label{EE:nstar} \notag
      \n^*\bet \ = \ (-1)^{j}*^{-1}\n*\bet \ + \
        (-1)^j*^{-1}\frac{d\rho}{\rho}\wedge*\bet.
\end{equation}
\end{lem}
\begin{proof}
Let us fix $\alp\in \Ome^{j-1}_{\comp}(M,E)$. Using \refe{Herm} we obtain
\begin{multline}\notag
      (\n\alp,\bet) \ = \ \int_M\, \<\n\alp\wedge*\obet\>\, \rho
      \\ = \
          \int_M\, d\<\alp\rho\wedge*\obet\>
      \ - \ (-1)^{j-1}\int_M\, \<\alp\wedge d\rho\wedge*\obet\>
      \ - \ (-1)^{j-1}\int_M\, \<\alp\rho\wedge\overline{\n*\bet}\>
      \\ = \
      (-1)^j\int_M\,
       \Big\<\, \alp\wedge*\,
        *^{-1}\big(\, \frac{d\rho}\rho\wedge*\obet\, \big)\, \Big\>\, \rho
      \ + \ (-1)^j\int_M\, \<\alp\wedge*\, *^{-1}\overline{\n*\bet}\>\rho
      \\ = \
      (-1)^j\Big(\, \alp,*^{-1}\n*\bet +*^{-1}\frac{d\rho}{\rho}\wedge*\bet\, \Big),
\end{multline}
where in the third equality we used that $ \int_M\, d\<\alp\rho*\obet\>= 0$ by the Stokes
theorem.
\end{proof}

To simplify the notation let us denote $A_\rho\bet= *^{-1}\frac{d\rho}{\rho}\wedge*\bet$.
Then
\begin{equation}\label{E:nstar2}
      \n^* \ = \ (-1)^j*^{-1}\n* + (-1)^jA_\rho:
      \Ome^j(M,E) \ \to \ \Ome^{j-1}(M,E).
\end{equation}
{}For the special case when $E=M\times\CC$ is the trivial line bundle, we obtain
\begin{equation}\label{E:dstar}
      d^* \ = \ (-1)^j*^{-1}d*+ (-1)^jA_\rho:\Ome^j(M) \ \to \ \Ome^{j-1}(M).
\end{equation}
In this paper we will use \refe{nstar2} and \refe{dstar} only in the case when $j=1$.
Clearly, for $\alp\in \Ome^0(M,E), \bet\in \Ome^j(M,E), \phi\in C^\infty(M)$ we have
\begin{equation}\label{E:A<>}
      A_\rho\<\alp\obet\> \ = \ \<\alp A_\rho\obet\>,
      \qquad A_\rho(\phi\bet) \ = \ \phi A_\rho\bet.
\end{equation}

\begin{defe}\label{D:Bochner}
The {\em Bochner Laplacian} is the operator $\n^*\n:C^\infty(E)\to C^\infty(E)$.
Similarly, we define the Laplacian $\Del_M=d^*d:C^\infty(M)\to C^\infty(M)$.
\end{defe}
Note, that if $d\mu$ is the Riemannian volume form on $M$, then $\Del_M= -\Del_g$, where
$\Del_g$ is the {\em metric} Laplacian, $\Del_gu=\div(\grad{}u)$.

\begin{rem}\label{R:g=g} In Sect.\ref{SS:Fmetric} we introduced a metric $g^{TM}$
associated with a first order elliptic differential operator $D$. It is easy to check that for
$D=\n:C^\infty(E)\to \Omega^1(M,E)$ this metric coincides with the metric $g$ which we now
assume to be given a priori.
\end{rem}

\subsection{Kato inequality}\label{SS:kato}
Recall that a distribution $\nu$ on $M$ is called positive (notation $\nu\ge0$), if for
every non-negative function $\phi\in C^\infty_c(M)$, we have $(\nu,\phi)\ge0$. It follows
that $\nu$ is in fact a positive Radon measure (see e.g. \cite{Gelfand-Vilenkin}, Theorem
1 in Sect. 2, Ch.II). We write $\nu_1\ge \nu_2$ if $\nu_1-\nu_2\ge0$.

The main result of this section is the following

\begin{thm}\label{T:kato}
Assume that $u\in\loloc(E)$ and $\n^*\n u\in\loloc(E)$.  Then
\begin{equation}\label{E:kato}
      \Del_M |u| \ \leq \ \RE\< \n^*\n u,\sign u\>,
\end{equation}
where
\[
      \sign u(x) \ = \
      \begin{cases}
          \frac{u(x)}{|u(x)|} & \textrm{if\ \ $u(x)\neq0$ },\\
          0 &\textrm{ otherwise}.
      \end{cases}
\]
\end{thm}

Let us emphasize that our Laplacians are {\em positive} operators. That explains the
discrepancy between \refe{kato} and the standard form of Kato inequality for scalar
functions on $\RR^n$, cf., for example, Theorem~X.27 of \cite{rs}.

The rest of this section is dedicated to the proof of \reft{kato}.

\subsection{Kato inequality in $C^{\infty}$ setting}\label{SS:smoothkato}
Let $u\in\esm$.  For $\eps>0$, set
\begin{equation}\label{E:epsabsvalue}
      |u|_{\eps} \ = \ \left(|u|^2+\epsilon^2\right)^{1/2}.
\end{equation}

\begin{Prop}\label{P:epsmoothkato}
Assume that $u\in\esm$.  Then the following inequality holds:
\begin{equation}\label{E:epsmoothkato}
      \uep\lapg\uep \ \leq \ \RE\< \n^*\n u,u\>.
\end{equation}
\end{Prop}
\begin{proof}
Let us fix $u\in\esm$. Using \refe{Herm}, we obtain
\begin{multline}\label{du2}
  2|u|_\eps d|u|_\eps \ = \ d|u|_\eps^2 \ = \ d |u|^2 \ = \ d\<u, u\>
  \ = \<(\nabla u)\bar u\>+\<u\left(\overline{\nabla u}\right)\>
  \\ = \
  2\RE \<(\nabla u)\bar u\> \ = \ 2\RE \<u\overline{\nabla u}\>,
\end{multline}
and hence
\begin{equation}\label{EE:du2a} \notag
|u|_\eps |d|u|_\eps| \ \le \ |u||\n u|.
\end{equation}
Since $|u|_\eps\ge |u|$, the above inequality implies
\begin{equation}\label{E:d<d}
  |d|u|_\eps|  \ \le \ |\n u|.
\end{equation}

{}Furthermore, using \refe{A<>}, we get
\begin{multline}\label{E:d*d}
      d^*(|u|_\eps d|u|_\eps) \ = \
       -*^{-1}d*(|u|_\eps d|u|_\eps) - |u|_\eps A_\rho d|u|_\eps
      \ = \ -*^{-1}d(|u|_\eps*d|u|_\eps) - |u|_\eps A_\rho d|u|_\eps
      \\ = \
       -*^{-1}\big(\, d|u|_\eps\wedge*d|u|_\eps\, \big)
                  -|u|_\eps*^{-1}d*d|u|_\eps - |u|_\eps A_\rho d|u|_\eps
       \\ = \ -\big|\, d|u|_\eps\, \big|^2 -|u|_\eps (*^{-1}d*+A_\rho)d|u|_\eps
       \ = \ -\big|\, d|u|_\eps\, \big|^2 +|u|_\eps \Del_M|u|_\eps.
\end{multline}
Similarly,
\begin{multline}\label{E:n*n}
      d^*\<u\overline{\n u}\> \ = \
       -*^{-1}d*\<u\overline{\n u}\> - A_\rho\<u\overline{\n u}\>
      \ = \
      -*^{-1}d\<u*\overline{\n u}\> - \<u A_\rho\overline{\n u}\>
      \\ = \
      -*^{-1}\<\n u\wedge*\overline{\n u}\> -*^{-1}\<u\overline{\n*\n u}\> -
  \<u A_\rho\overline{\n u}\>
      \ = \ - \big|\, \n u\, \big|^2+ \<u\overline{\n^*\n u}\>.
\end{multline}
{}From \eqref{du2}, \refe{d*d} and \refe{n*n} we obtain the equality
\[
       -\big|\, d|u|_\eps\, \big|^2 +|u|_\eps \Del_M|u|_\eps
       \ = \
       - \big|\, \n u\, \big|^2+ \RE\<u\overline{\n^*\n u}\>
=- \big|\, \n u\, \big|^2+ \RE\<(\n^*\n  u)\bar u\>,
\]
which, in view of \refe{d<d}, implies \refe{epsmoothkato}.
\end{proof}

\begin{cor}\label{C:smoothkato}
Assume that $u\in\esm$.  Then the following distributional inequality holds:
\begin{equation}\nonumber
      \lapg |u| \ \leq  \ \RE\< \n^*\n u,\znak u\>.
\end{equation}
\end{cor}
\begin{rem}
In the case when $d\mu$ is the Riemannian volume form on $M$ (so that $\rho\equiv1$)
analogues of \refp{epsmoothkato} and \refc{smoothkato} were obtained by H.~Hess,
R.~Schrader and D.~A.~Uhlenbrock \cite{hsu2}. They used a slightly more complicated
definition of $\sign{u}$, which is, in fact, unnecessary, because the difference with our
definition occurs only on a set of measure 0.
\end{rem}

\subsection{Friedrichs mollifiers}\label{SS:friedrichs}
We will now deduce \reft{kato} from \refp{epsmoothkato} using the Friedrichs mollifiers
technique.  For reader's convenience we will review basic definitions and results of
K.~Friedrichs~\cite{fri}.

Suppose that $j\in C_c^{\infty}(\RR^n)$, $j(z)\ge 0$ for all $z\in\RR^n$, $j(z)=0$ for
$|z|\geq1$, and $\int_{\RR^n}j(z)\,dz=1$.

{}For $\rho>0$ and $x\in\RR^n$, define $j_\rho(x)=\rho^{-n}j(\rho^{-1}x)$. Then
$j_\rho\in\rncomp$, $j_\rho\ge 0$, $j_{\rho}(x)=0$  if $|x|\ge \rho$, and $\int_{\RR^n}
j_{\rho}(x)\,dx=1$.

Let $\Omega\subset\RR^n$ be an open set.  Suppose that $f\in\loloc(\Omega)$. Define
\begin{equation}\label{E:avg}
    (J^{\rho}f)(x)  \ = \  f^{\rho}(x) \ = \  \int\, f(x-y)j_{\rho}(y)\,dy
     \ = \  \int\, j_{\rho}(x-y)f(y)\,dy,
\end{equation}
where the integration is taken over $\RR^n$. Then $J^{\rho}f(x)$ is defined for
$x\in\Omega_{\rho}$, where
 $$
   \Omega_{\rho} \ = \ \big\{\, x:\,
     x\in\Omega,\dist(x,\partial\Omega) > \,\rho\, \big\}.
$$
 We recall the following standard Lemma (cf. \cite[Theorem 1.3.2]{Horm1}):
%
\begin{lem}\label{L:mollifier}
Let $1\le p<\infty$.  With the notations from \refe{avg},
\begin{enumerate}
\item If $f\in\lploc(\Omega)$, then $J^{\rho}f\in C^{\infty}(\Omega_\rho)$
\item If $f\in\lploc(\Omega)$, then $J^{\rho}f\to f$ as $\rho\to 0$,
in the norm of $L^{p}(K)$, where $K\subset\Omega$ is any compact set.
\item If $f\in C(\Omega)$, then $J^{\rho}f\to f$ as $\rho\to 0$,
uniformly on any compact set $K\subset\Omega$.
\end{enumerate}
\end{lem}

To apply the Friedrichs mollifiers technique to the operator $\n^*\n$, let us choose a
coordinate neighborhood $U$ of $M$ with coordinates $x^1,x^2,\dots x^{n}$. Choosing a
local frame of $E$ over $U$ we will identify the space of sections of $E$ over $U$ with
the space of vector valued functions on $U$. Under this identification the Bochner
Laplacian $\n^*\n$ becomes a second order elliptic differential operator
which can be written in the form
\begin{equation}\label{EE:katopaper}\notag
      \sum_{i,k}a_{ik}(x)\pari\park+\sum_{i}b_i(x)\pari+c(x),
\end{equation}
where $a_{ik}$, $b_i$ and $c$ are $m\times m$ matrices of smooth functions ($m=\dim E$).
Alternatively, we can write it in the form
\begin{equation}\label{E:katopaper}
      \sum_{i,k}\pari\, a_{ik}(x)\park+\sum_{i}b_i(x)\pari+c(x),
\end{equation}
with possibly different $a_{ik}$, $b_i$ and $c$. The later form is sometimes more
convenient if the coefficients are not smooth.

We can also consider $u^{\rho}= \mathcal{J}^{\rho}u$, where $\sj$ is an integral operator
whose integral kernel is $\kjr(x-y)\Id$, where $\Id$ is $m\times m$ identity matrix.

We now state a key proposition whose proof is given in Appendix~A.

\begin{Prop}\label{P:keykato}
Assume that $u\in\loloc(E)$ and $\n^*\n{}u\in\loloc(E)$. Then $\n^*\n \ur\to \n^*\n{}u$
in $\loloc(E)$ over $U$, as $\rho\to 0+$.
\end{Prop}

\subsection{Proof of \reft{kato}}\label{SS:prkato}
Note first that the statement is actually local. Namely,  using partition of unity we see
that it suffices to prove the inequality over every coordinate neighborhood from a
covering of $M$. Let us fix such a neighborhood $U$ and apply \refp{keykato}. Clearly
$\ur\in C^\infty(E|_{U_\rho})$.
  By \refp{epsmoothkato}, for any $\epsilon >0$,
\begin{equation}\label{E:prelimkato}
      \lapg |\ur|_{\epsilon} \ \le \ \RE \Big\<
      \n^*\n\ur,\frac{\ur}{|\ur|_{\eps}}\Big\>.
\end{equation}
In this inequality we will first fix $\eps>0$ and
pass to the limit as $\rho\to 0+$.

Consider the left hand side of \refe{prelimkato}. Clearly,
\begin{equation}\label{E:pointkato}
      \big|\, |\ur|_{\epsilon}-|u|_{\epsilon}\, \big| \ \leq \
      \big|\, |\ur|-|u|\, \big| \ \leq \ |\ur-u|.
\end{equation}
Since $\ur\to u$ in $\loloc(E|_{U})$, it follows from \refe{pointkato} that
$|\ur|_\epsilon\to |u|_{\epsilon}$  in $\loloc(U)$ as $\rho\to 0+$.  This immediately gives
$\lapg |\ur|_\epsilon\to \lapg |u|_\epsilon$ in distributional sense over $U$ (or, more
precisely, over any relatively compact open subset).
Now we turn to the right hand side of \refe{prelimkato}.  Our goal is to show that
\begin{equation}\label{E:rhkato}
      \RE\Big\<\, \n^*\n\ur,\qrkato\, \Big\> -
       \RE\Big\<\, \n^*\n u,\qkato\, \Big\>  \to \  0
       \qquad\text{in}\quad \loloc(U),
\end{equation}
as $\rho\to 0+$.
By adding and subtracting the term $\RE\<\n^*\n u,|u^\rho|_\eps^{-1}u^\rho\>$ at the
left hand side of~\refe{rhkato}, we get
\begin{equation}\label{E:katosumrule}
      \RE\Big\<\,  \n^*\n \ur- \n^*\n u,\qrkato\, \Big\> \ + \
      \RE\Big\<\, \n^*\n u,\qrkato-\qkato\, \Big\>.
\end{equation}

Let $\rho\to 0+$.  By \refp{keykato}, $\n^*\n \ur\to \n^*\n u$ in
$\loloc(E|_U)$.  We also know that $\left||u^\rho|_\eps^{-1}u^\rho\right|<1.$
Therefore, the first term of the sum in \refe{katosumrule} converges to $0$ in
$\loloc(U)$ (hence in distributional sense).

We know that $\n^*\n u\in\loloc(E)$.  We have already shown that as $\rho\to 0+$,
$|\ur|_{\epsilon}\to|u|_{\epsilon}$ in $\loloc(U)$, and hence, after passing to a
subsequence, almost everywhere. Therefore, along the same subsequence,
$|u^\rho|_{\epsilon}^{-1}u^\rho \to|u|_\epsilon^{-1}u$ almost everywhere. So the
dominated convergence theorem implies that the second term of the sum in
\refe{katosumrule} goes to $0$ in $\loloc(U)$ (hence in distributional sense ) along a sequence
of positive $\rho$'s.

Therefore, for all $u\in\loloc(E)$ such that $\n^*\n  u\in\loloc(E)$, we obtain
\begin{equation}\label{E:katoepsilon}
      \lapg |u|_\epsilon \ \leq  \ \Realna \Big\<\,  \n^*\n u,\qkato\, \Big\>,
\end{equation}

In this inequality we will pass to the distributional limit as $\epsilon\to 0+$.
Let us take a
real-valued $\phi\in\mcomp$, and consider
\begin{equation}\label{E:repskato}
    \int\, \phi \lapg |u|_\epsilon \,d\mu
    \ = \
    \int\, |u|_{\epsilon}\lapg\phi\, d\mu \ \to \ \int\, |u|\lapg\phi\,d\mu.
\end{equation}
The limit in \refe{repskato} can be understood as the value of the distribution
$\lapg |u|$ on a test function~$\phi$.
Therefore, $\lapg |u|_\epsilon\to\lapg |u|$ in distributional sense,
as $\epsilon\to 0+$.

On the right hand side of \refe{katoepsilon}, $|u|_\epsilon^{-1}u\to \znak u$
almost everywhere and with a uniform bound $\left||u|_\epsilon^{-1}u\right|\leq 1$.
Since $\n^*\n u\in\loloc(E)$, the dominated convergence theorem shows that
the right hand side converges to $\RE\big<\n^*\n  u,\znak u\big>$ in
$\loloc(M)$ (and hence in distributional sense).~$\square$

\section{The case when the support of $V$ is small}\label{S:smallsup}

In this section we prove a particular case of \reft{main}, when $D=\n$ and the potential
$V$ is supported in a coordinate neighborhood $W\subset M$.

\begin{propos}\label{P:smallsup}
Suppose $M$ is a manifold endowed with a Riemannian metric $g^{TM}$ and a smooth measure
$d\mu$. Let $W\subset M$ be a relatively compact coordinate neighborhood in $M$ and
suppose that $V=V_++V_-$ is as in the Assumption~A. Suppose, in addition, that
$\supp{V}\subset W$. Then the operator
\[
      H_V \ = \ \n^*\n+V
\]
is essentially self-adjoint on $C^\infty_c(E)$.
\end{propos}
We precede the proof with several lemmas, which will be also used in the subsequent
sections.
\begin{lem}\label{L:hV-<}
If the Assumption~A is satisfied and $\supp{V_-}\subset W$, then there exist positive
constants $a<1, \ C_a$ such that
\begin{equation}\label{E:hV-<}
      |(V_-u,u)| \ \le \ a\, \|\n u\|^2 + C_a\|u\|^2,
      \qquad u\in C^\infty_c(E).
\end{equation}
\end{lem}
\begin{proof}
Let $\oW$ denote the closure of $W$ in $M$.
By \refe{V<Del}, there exist $a_\oW, \ 0\le a_{\oW}<1$, and $C_\oW>0$,
such that
\begin{equation*}
      \big\||V_-|v\big\| \ \le \
       a_\oW\, \big\|\, \Del_M v\, \big\|+
            C_\oW\, \big\|\, v\, \big\|,
      \quad v\in C^\infty_c(M).
\end{equation*}
In other words the operator $|V_-|$ is dominated by $\Del_M$. Hence by Theorem X.18
from \cite{rs} the quadratic form of $|V_-|$ will be dominated by the quadratic form of $\Del_M$.
More precisely, for any $a\in (a_\oW, 1)$
there exists $C_a>0$ such that
\begin{equation*}
      \big(|V_-|v,v\big) \ \le \
       a\big(\Del_M v,v\big)+C_a\|v\|^2=a\|dv\|^2+C_a\|v\|^2, \quad  v\in C^\infty_c(M).
\end{equation*}
Let us apply this to $v=|u|_\eps$, where $u\in C^\infty_c(E)$ and $|u|_\eps$
is defined as in \refe{epsabsvalue}. We obtain then, using the inequality \refe{d<d}:
\begin{equation*}
|(V_- u, u)|\le  (|V_-||u|_\eps,|u|_\eps)
\le a\|d|u|_\eps\|^2+C_a\big\||u|_\eps\big\|^2\le a\|\n u\|^2+C_a \big\||u|_\eps\big\|^2.
\end{equation*}

Taking the limit as $\eps\to 0+$, we come to \refe{hV-<}.
\end{proof}
\begin{cor}\label{C:boundedbelow}
If the Assumption~A is satisfied and $\supp{V_-}\subset W$, then the operator $H_V$ is
semi-bounded below on $C^\infty_c(E)$, i.e., there exists $C>0$ such that
  \begin{equation}\label{EE:boundedbelow} \notag
     (H_Vu,u) \ \ge \ -C\|u\|^2, \qquad \text{for all}\quad u\in C^\infty_c(E).
  \end{equation}
\end{cor}

The following lemma follows immediately from \refe{defsym}.
\begin{lem}\label{L:commut}
Suppose $H_V=D^*D+V$ is as in \refe{HV}. For every $\psi\in C^2(M)$ and every $u\in
L^2_{\loc}(E)$ the following distributional equality holds
\begin{equation}\label{EE:commut} \notag
      H_V(\psi u) \ = \
      \psi H_Vu \ + \ D^*\big(\, \hatD(d\psi)u\, \big)
       \ - \ \big(\, \hatD(d\psi)\, \big)^*Du.
\end{equation}
\end{lem}
\begin{cor}\label{C:commut}
Suppose $u\in \Dom(H_{V,\max})$. Let $\psi\in C^2(M)$ be such that $\supp(d\psi)$
is compact. If the restriction of $u$ to an open neighborhood of $\supp(d\psi)$ is in
$W^{1,2}_{\loc}$, then $\psi{u}\in \Dom(H_{V,\max})$.
\end{cor}

\subsection{Proof of \refp{smallsup}}\label{SS:prsmallsup}
{}Fix an open neighborhood $U\subset W$ of $\supp{V}$ such that the closure $\oU$ of $U$
is contained in $W$. Let $\phi,\tilphi:M\to [0,1]$ be smooth functions such that
$\phi^2+\tilphi^2\equiv1$, the restriction of $\phi$ to $U$ is identically equal to 1,
and $\supp\phi\subset W$. It follows that $\supp\tilphi\subset M\backslash{U}$.

By IMS localization formula (cf. \cite[\S3.1]{CFKS}, \cite[Lemma~3.1]{sh4}), we have
\begin{equation}\label{E:IMS1}
     \n^*\n
      \ = \
      \phi\n^*\n\phi+\tilphi\n^*\n\tilphi+
       \frac{1}{2}\left[[\n^*\n,\phi],\phi\right]+
          \frac{1}{2}\left[[\n^*\n,\tilphi],\tilphi\right],
\end{equation}
where $[\cdot,\cdot]$ denotes the commutator bracket. Set
\[
      J  \  = \ \frac{1}{2}\left[[\n^*\n,\phi],\phi\right]+
          \frac{1}{2}\left[[\n^*\n,\tilphi],\tilphi\right]
      \ = \
      -\frac12|d\phi|^2-\frac12|d\tilphi|^2.
\]
Then $J$ is a smooth function on $M$, whose support is contained in $W\backslash{U}$.
With this notation \refe{IMS1} takes the form
\begin{equation}\label{E:IMS3}
     \n^*\n
      \ = \
      \phi\n^*\n\phi+\tilphi\n^*\n\tilphi+J.
\end{equation}

Let $C>0$ be as in \refc{boundedbelow}. Let $b\gg C$ be a large number, which will be
chosen later. Suppose that $u\in L^2(E)$ satisfies the equality
\begin{equation}\label{E:H+b=0}
      (H_V+b)u \ = \ 0.
\end{equation}
Here the expression $H_Vu$ is understood in the sense of distributions. To prove the
theorem, it is enough to show that $u=0$ (cf., e.g., Theorem~X.26 of \cite{rs} or
\cite{Glazman}). Arguing {\em ad absurdum}, let us assume that $u\not=0$.

Since the restriction $V|_{M\backslash{U}}$ of $V$ to $M\backslash{U}$ vanishes,
\refe{H+b=0} implies that $\n^*\n{u}|_{M\backslash{U}}= -bu|_{M\backslash{U}}\in L^2(E)$.
Hence, $u|_{M\backslash{U}}\in W^{2,2}_{\loc}$ by elliptic regularity. By \refc{commut},
the sections $\phi{u}$ and $\tilphi{u}$ belong to $\Dom(H_{V,\max})$. Note also that on
the support of $\tilphi$ we have $H_V= \n^*\n$. It is well known (cf., for example,
\cite{hsu2}), and also follows from \refp{pospot} below, that the operator $\n^*\n$ is
essentially self-adjoint (note that Assumption~A of \refss{setting} is tautologically
true in this case, while Assumption~B follows from elliptic regularity). Hence,
by taking closure (see also arguments in the proof of Proposition \ref{P:posquadform})
we obtain
\begin{equation}\label{E:tilphi>}
      (\tilphi\n^*\n(\tilphi u), u) \ = \ \|\n(\tilphi u)\|^2\ge 0
\end{equation}
and the equality is attained only if $\tilphi{u}=0$.

Let $\tilb=b-\max_{x\in M}|J(x)|$ and assume that $b$ is large enough so that $\tilb>0$.
>From \refe{IMS3} and \refe{H+b=0}, we obtain
\begin{equation}\label{E:noJ}
    0 \ = \ \Big(\, (H_V+b)u,u\, \Big) \ \ge \
    \Big(\, \phi H_V(\phi u),u\, \Big)  +
    \Big(\, \tilphi\n^*\n(\tilphi u), u\, \Big) + \tilb\, \|u\|^2.
\end{equation}
Since we assumed $u\not=0$, the equations  \refe{tilphi>} and \refe{noJ} imply that
\begin{equation}\label{E:H+b<0}
      \Big(\, \phi(H_V+\tilb)(\phi u),u\, \Big) \ \le  \
       \Big(\, \phi H_V(\phi u),u\, \Big) \ + \ \tilb\, \|u\|^2 \ \le \ 0,
       \qquad \phi u\not=0.
\end{equation}

We can and we will view the coordinate neighborhood $W$ as a subset of $\RR^n$. Fix a
Riemannian metric $\tilg$ on $\RR^n$, whose restriction to $W$ coincides with the metric
induced by $g^{TM}$, and which is equal to the standard Euclidean metric on $\RR^n$
outside of some compact set $K\supset{W}$. We also endow $\RR^n$ with a measure $\tilmu$,
whose restriction to $W$ coincides with the measure induced by the measure $\mu$ on $M$,
and whose restriction to $\RR^n\backslash{K}$ is the Lebesgue measure on $\RR^n$.

Let $\tilE:=\RR^n\times\CC^m$ be the trivial vector bundle over $\RR^n$ endowed with the
standard Hermitian metric. Choosing an orthonormal frame of the restriction $E|_W$ of $E$
to $W$ we can and we will identify $E|_W$ and $\tilE|_W$. Then the connection $\n$
induces a Hermitian connection on $\tilE|_W$. Let $\tiln$ be a Hermitian connection on
$\tilE$, whose restriction to $W$ coincides with the connection induced by $\n$, and
whose restriction to $M\backslash{K}$ is the trivial flat connection. Then, by
\refe{H+b<0}, we have
\[
      \Big(\, \phi(\tiln^*\tiln+V+\tilb)(\phi u),u\, \Big) \ \le \ 0,
      \qquad \phi u\not=0.
\]
Let us choose $\tilb$ large enough so that $\tilH_V= \tiln^*\tiln+V+\tilb>0$ on
$C^\infty_c(\tilE)$ (this is possible due to \refc{boundedbelow}). It follows (cf.
Theorem~X.26 of \cite{rs}) that the operator $\tilH_V$ {\em is not} essentially
self-adjoint on $C^\infty_c(\tilE)$. Thus we reduced the proof of the proposition to the
analogous question for the operator $\tilH_V$ on $\RR^n$.

Since the operator $\tilH_V$ is not essentially self-adjoint, we conclude from
Theorem~X.26 of \cite{rs} that there exists non-zero $v\in L^2(\tilE)$ such that
\begin{equation}\label{E:tilH=0}
      (\tilH_V+\tilb)\,v \ = \ 0.
\end{equation}
Let $\tilDel:=d^*d$ denote the scalar Laplacian on $\RR^n$ associated to the metric
$\tilg$ and the measure $\tilmu$, cf. \refd{Bochner}. From \refe{tilH=0} and the Kato
inequality \refe{kato} we obtain
\[
      \tilDel|v| \le \RE\Big\<\, \tiln^*\tiln v,\sign v\, \Big\>
      \ = \
      -\Big\<\, (V+\tilb)v,\sign v\, \Big\> \ \le \   (|V_-|-\tilb)\, |v|,
\]
where $|V_-(x)|$ is the norm of the endomorphism $V_-(x):\tilE_x\to \tilE_x$. Thus
\begin{equation}\label{E:del<}
      (\tilDel+\tilb)|v| \ \le \ |V_-|\, |v|.
\end{equation}
It is well known and also explained in Appendix~B that the Schwartz kernel of the
operator $(\tilDel+\tilb)^{-1}$ is positive. A somewhat more subtle argument, cf.
\refp{bounded-geometry} of Appendix~B, shows that the distributional inequality
\refe{del<} implies
\begin{equation}\label{E:u<}
      |v| \ \le (\tilDel+\tilb)^{-1}|V_-|\, |v|,
\end{equation}
where the inverse operator $(\tilDel+\tilb)^{-1}$ is understood , e.g., as an operator on
the space $\mathcal{S}'(\RR^n)$ of tempered distributions.

Let $\oW$ denote the closure of $W$ in $M$. Since $\supp{}V_-\subset \oW$, it follows
from \refe{V<Del}, that
\[
      \big\|\, |V_-|s\, \big\| \ \le \ a_\oW\|\tilDel s\|+ C_\oW\|s\|,
      \qquad\text{for all}\quad s\in C^\infty_c(\RR^n).
\]
Hence,
\begin{equation*}\label{EE:Vdel} \notag
      \big\|\, |V_-|(\tilDel+\tilb)^{-1}s\, \big\|
     \ \le \
     a_\oW\, \big\|\, \tilDel(\tilDel+\tilb)^{-1}s\, \big\| +
     C_\oW\, \big\|\, (\tilDel+\tilb)^{-1}s\, \big\|
     \ \le \
     \big(\, a_{\oW}+C_\oW/\tilb\, \big)\|s\|.
\end{equation*}
Thus we can choose $\tilb$ large enough so that $\big\||V_-|(\tilDel+\tilb)^{-1}\big\|<
1$. Since
\[
      \big(\, (\tilDel+\tilb)^{-1}|V_-|\, \big)^*
      \ = \ |V_-|(\tilDel+\tilb)^{-1},
\]
we conclude that $\|(\tilDel+\tilb)^{-1}|V_-|\|<1$. Then \refe{u<} together with our
assumption $v\not=0$ implies  $\|v\|<\|v\|$. The obtained contradiction proves the
proposition. \hfill$\square$

%
\begin{rem}\label{R:simple}
One of the main steps of the above proof was the reduction to the case $M=\RR^n$, where
we can use \refp{bounded-geometry} to conclude that the inequality \refe{del<} implies
\refe{u<}. The proof would be much simpler if we knew that the same is true on arbitrary
complete Riemannian manifold, i.e., if we knew that Conjecture~P of Appendix~B holds.
Unfortunately, we don't know if this is the case (see Appendix~B, where we explain the
difficulties one meets in an attempt to prove this conjecture).
\end{rem}


\refp{smallsup} and \refc{boundedbelow} imply, cf. Theorem~X.26 of \cite{rs}, the
following
\begin{cor}\label{C:bounded}
Let $V, \ C$ be as in \refc{boundedbelow}. Then
\[
     (H_Vu,u) \ \ge \ - C\|u\|^2, \qquad \text{for all} \quad u\in
\Dom(H_{V,\max}).
\]
\end{cor}

Let $W^{1,2}(E)$ denote the completion of the space $C^\infty_c(E)$ with respect to the
norm $\|\cdot\|_{{}_{W^{1,2}}}$ defined by the scalar product
\begin{equation}\label{E:W12norm}
    (u,v)_{{}_{W^{1,2}}} \ := \ (u,v) \ + \ (\n u,\n v)
    \qquad u, v\in C^{\infty}_c(E).
\end{equation}
(It is not difficult to see that $W^{1,2}(E)$ coincides with the set of all $u\in
L^2(E)$, such that $\n{}u\in L^2(T^*M\otimes{}E)$, but we will not use this fact.)

Since the operator $H_V$ is essentially self-adjoint,
the domain $Q(H_V)$ of its quadratic form coincides
with the closure of the space $C^\infty_c(E)$ with respect to the norm $\|\cdot\|_1$
defined by the formula
\[
      \|u\|_1^2 \ := \ (H_Vu,u)+(C+1)\|u\|^2.
\]
By Theorem~X.23 of \cite{rs}, we have $\Dom(H_{V,\max})\subset Q(H_V)$.
\begin{cor}\label{C:domH}
Under the assumptions of \refp{smallsup}, we have $Q(H_V)\subset W^{1,2}(E)$. In
particular, $\Dom(H_{V,\max})\subset W^{1,2}(E)$.
\end{cor}
\begin{proof}
By \refl{hV-<}, there exist $a<1$ and $C_a>0$ such that
\[
      (H_Vu,u) \ = \  \|\n u\|^2+(V_+u,u) + (V_-u,u) \ \ge \
      (1-a)\|\n u\|^2-C_a\|u\|^2,
     \qquad u\in
C^\infty_c(E).
\]
Thus the domain $Q(H_V)$ of the quadratic form of the operator $H_V$ coincides with the
closure of the space $C^\infty_c(E)$ with respect to the norm
\[
     \|u\|_1^2 \ : = \ (H_Vu,u) + (C_a+1)\|u\|^2 \ \ge \ (1-a)\|\n u\|^2+\|u\|^2.
\]
Hence $Q(H_V)$ is contained in the closure of $C^\infty_c(E)$ with respect to the norm
$(1-a)\|\n u\|^2+\|u\|^2$, which coincides with the space $W^{1,2}(E)$.
\end{proof}

\section{Proof of the second part of \reft{domain}}\label{S:domain}

In this section we assume that there exists a $C^\infty$ function $s:T^*M\to \RR$ such
that
\[
     \sig(D^*D)(\xi) \ = \ s(\xi)\Id,
\]
for each $\xi\in T^*M$. Then,
\[
     \sig(D^*D) \ = \ |\xi|^2 \Id,
\]
where $|\xi|$ denotes the norm on $T^*M$ defined by the scalar product \refe{trace}. In
this situation we say that the operator $D^*D$ {\em has a scalar leading symbol}. It
follows from Proposition~2.5 of \cite{bgv}, that there exists a Hermitian connection $\n$
on $E$ and a linear self-adjoint bundle map $F\in C^\infty(\End{E})$ such that
\begin{equation}\label{EE:DtD=nn} \notag
     D^*D \ = \ \n^*\n + F.
\end{equation}

\begin{propos}\label{P:domain}
Let $D$ be as above and consider the operator
\[
     H_V \ = \ D^*D+V \ = \ \n^*\n + F+V,
\]
where $V=V^++V^-$ be as in Assumption~A. Then $\Dom(H_{V,\max})\subset\wl$.
\end{propos}
\begin{proof}
Let $u\in \Dom(H_{V,\max})$. We need to prove that $u\in\wl$. Clearly, it is enough to
show that for every $x\in M$ there is an open coordinate neighborhood $U\ni x$ and a
constant $C_U$ such that for all $j=1\nek n$ and all $z\in C^\infty_c(E|_U)$ we have
\begin{equation}\label{EE:unjz} \notag
     \big|\, (\, \n_{\frac{\p}{\p x_j}}u,\, z\, )\, \big|
     \ \le \ C_U\, \|z\|.
\end{equation}
Let us fix $x\in M$ and choose a coordinate neighborhood $W_1$ of $x$. Let $\phi:M\to
[0,1]$ be a smooth function, whose support is contained in $W_1$, which is identically
equal to 1 near $x$. In other words we assume that there exists an open neighborhood
$W_2$ of $x$, such that $\phi|_{W_2}\equiv1$. Set
\[
     H_\phi \ := \ \n^*\n+\phi (F+V)\phi.
\]
Since $F$ is smooth, the potential $\phi(F+V)\phi$ satisfies Assumption~A. Hence, by
\refp{smallsup}  and Corollaries~\ref{C:bounded} and \ref{C:domH}, the operator $H_\phi$
is essentially self-adjoint, bounded below, and $\Dom(H_{\phi, \max})\subset W^{1,2}(E)$.

Choose an open neighborhood $U\subset W_2$ of $x$ and a  smooth function $\psi:M\to
[0,1]$ whose support is contained in $W_2$, such that $\psi|_U\equiv1$.
\begin{lem}\label{L:Hp<}
There exist a constant $C_1=C_1(u,\phi,\psi)$ such that
\begin{equation}\label{EE:Hp<} \notag
     |(\psi u,H_\phi s)| \ \le \ C_1\, \big(\,
         \| s\|+\|\n s\|\, \big),
    \qquad\text{for all}\quad s\in \Dom(H_{\phi, \max}).
\end{equation}
\end{lem}
\begin{proof}
{}Fix $s\in \Dom(H_{\phi, \max})$. Since $H_\phi$ is essentially self-adjoint on
$C^\infty_c(E)$, there exists a sequence $s_k\in C_c^\infty(E)$, which converges to $s$
in the graph-norm of $H_\phi$. Since $\Dom(H_{\phi, \max})\subset W^{1,2}(E)$, we have
$s_k\to s$ also in the topology of $W^{1,2}(E)$. Hence,
\[
     \lim_{k\to\infty}\, |(\psi u, H_\phi s_k)| \ = \ |(\psi u, H_\phi s)|;
     \qquad
     \lim_{k\to\infty}\,
     \Big(\, \|s_k\|+\|\n s_k\|\, \Big) \ = \ \|s\|+\|\n s\|.
\]
It follows that it is enough to prove the lemma for the case when $s\in C^\infty_c(E)$,
which we will henceforth assume.

Since the support of $\psi$ is compact, we conclude from \refl{commut} that
$\psi{}H_V-H_V\psi$ defines a continuous map $W^{1,2}(E)\to L^2(E)$. Hence, there exists a
constant $C>0$ such that
\begin{equation}\label{E:pH-Hp}
     \|\psi H_Vs-H_V(\psi s)\|  \ \le \ C(\|s\|+\|\n s\|).
\end{equation}
Since $\phi|_{\supp \psi}\equiv1$ we have $\psi{}H_\phi=\psi{}H_V$. Thus, using
\refe{pH-Hp}, we obtain
\begin{multline}\notag
     |(\psi u, H_\phi s)| \ = \ |(u,\psi H_V s)| \ \le \
     |(u,H_V(\psi s))| + C(\| s\|+\|\n s\|)\|u\|
     \ \le \
     |(H_Vu,\psi s)| \\ + C(\| s\|+\|\n s\|)\|u\| \ \le \
     \|H_Vu\|\|\psi s\| + C(\|s\|+\|\n s\|)\|u\| \ \le \
     C_1 (\| s\|+\|\n s\|),
\end{multline}
where $C_1=C\|u\|+\|H_Vu\|$.
\end{proof}
\begin{lem}\label{L:Hp<2}
There exist positive constants $a<1, C_2$ such that
\begin{equation}\label{E:Hp<2}
     (H_\phi s,s) \ \ge \ (1-a)\|\n s\|^2-C_2\|s\|^2,
     \qquad   s\in \Dom(H_{\phi, \max}).
\end{equation}
\end{lem}
\begin{proof}
As in the proof of \refl{Hp<}, it is enough to prove \refe{Hp<2} for $s\in
C^\infty_c(E)$, which we will henceforth assume. By \refl{hV-<}, there exist positive
constants $a'<1, \ C'$ such that
\begin{equation}\label{E:V-<3}
     \left|\, \big(\, \phi V_-(\phi s),\, s\, \big)\, \right|
     \ \le \ a'\, \big\|\, \n(\phi s)\, \big\|^2+ C'\|s\|^2,
     \qquad  s\in C^\infty_c(E).
\end{equation}
{}Fix $\eps>0$ such that $a:=a'(1+\eps)<1$. Set $m_1:= \max_{x\in M}|d\phi(x)|$. Then
\begin{equation}\label{E:np<}
     \big\|\, \n(\phi s)\, \big\|^2 \ \le \
       \big(\, \|\n s\|+m_1\|s\|\, \big)^2 \ \le \
     (1+\eps)\, \|\n s\|^2+
            (1+1/\eps)\, m_1^2\, \|s\|^2.
\end{equation}
Combining \refe{V-<3} and \refe{np<} we obtain
\begin{equation}\label{E:V-<4}
     \left|\, \big(\, \phi V_-(\phi s),\, s\, \big)\, \right|
     \ \le \ a\, \|\n s\|^2+ C''\|s\|^2,
\end{equation}
where $C''=C'+a'(1+1/\eps)m_1^2$. Set $m_2=\max_{x\in M}|\phi{}F\phi|$. From \refe{V-<4}
and the inequality $V_+\ge0$, we obtain
\begin{multline}\notag
     (H_\phi s,s) \ \ge \ (\n^*\n s,s) - |(\phi V_-\phi s,s)|
              - |(\phi{F}\phi s,s)|
     \ \ge \ (1-a)\|\n s\|^2 - (C''+m_2)\|s\|^2.
\end{multline}
\end{proof}

We return to the proof of \refp{domain}. Fix $z\in C_c^\infty(E|_U)$ and
$j\in\{1\nek{}n\}$. Since the operator $H_\phi$ is essentially self-adjoint and bounded
below, for sufficiently large $\lam\gg0$, the map $H_\phi+\lam:\Dom(H_{\phi, \max})\to
L^2(E)$ is bijective. Set
\[
     s \ = \ (H_\phi+\lam)^{-1}\n_{\frac{\p}{\p x_j}}^*z,
\]
where $\n_{\frac{\p}{\p x_j}}^*$ denotes the formal adjoint of $\n_{\frac{\p}{\p x_j}}$.
Then, using \refl{Hp<2}, we obtain
\begin{equation}\label{E:Hpl>}
    \big(\, (H_\phi+\lam)s,s\, \big) \ \ge \  (1-a)\|\n s\|^2 + (\lam-C_2)\|s\|^2.
\end{equation}
Thus, if $\lam>C_2$, then $\big(\, (H_\phi+\lam)s,s\, \big)>0$.  Since $s\in
\Dom(H_{\phi,\max})\subset W^{1,2}(E)$, we get
\[
     \big(\, (H_\phi+\lam)s,s\, \big) \ = \
     (\n_{\frac{\p}{\p x_j}}^*z,s) \ = \ (z,\n_{\frac{\p}{\p x_j}}s).
\]
{}For every $\eps>0$, we have
\begin{multline}\label{E:2}
       \big(\, (H_\phi+\lam)s,s\, \big)
       \ = \ (z,\n_{\frac{\p}{\p x_j}}s)
        \ \le \  \|z\|\,
        \Big(\, \int_U\, |\n_{\frac{\p}{\p x_j}}s|^2\, d\mu\, \Big)^{1/2}
       \\ \le \
       \frac{\eps}2\, \int_U\, |\n_{\frac{\p}{\p x_j}}s|^2\, d\mu
       +\frac1{2\eps}\|z\|^2
        \ \le \
        \frac{\eps C_3}2\|\n s\|^2+\frac1{2\eps}\|z\|^2,
\end{multline}
where $C_3$ is the supremum of the length of the vector $\frac{\p}{\p x_j}$ on $U$.
Combining \refe{Hpl>} and \refe{2}, we get
\begin{equation}\label{E:3}
     \|\n s\| +\|s\| \ \le \ C_4\|z\|
\end{equation}
for some constant $C_4$. Using \refl{Hp<} and \refe{3} we get
\[
     |(u,\n_{\frac{\p}{\p x_j}}^*z)|
     \ = \ |(\psi u,\n_{\frac{\p}{\p x_j}}^*z)| \ = \
     |(\psi u, (H_\phi+\lam)s| \ \le \ C_1(\|\n s\|+\|s\|)+ C_5\|s\|
     \ \le \ C_6\|z\|,
\]
where $C_5=\|\psi{}u\|$ and $C_6= (C_1+C_5)C_4$.
\end{proof}

\section{Quadratic forms and the essential self-adjointness
                        of $H_{\vp}$}\label{S:quadratic}

Suppose that the operator $H_V=D^*D+V$ satisfies Assumptions~A and B.
Throughout this section we assume that
\begin{equation}\label{E:HV>}
     (H_Vu,u) \ \ge \ \|u\|^2, \qquad \text{for all}\quad
     u\in C^\infty_c(E).
\end{equation}

\subsection{}\Label{SS:formhV}
Recall that $W^{-1,2}_{\loc}(E)$ is the dual space to $W^{1,2}_{\comp}(E)$.
It follows that the intersection
$W^{-1,2}_{\loc}(E)\cap{}L^1_{\loc}(E)$
consists of the sections
$v\in L^1_{\loc}(E)$ satisfying the following property: for each compact set $K\subset
M$, there exists a constant $C_{K,v}>0$ such that
\[
    \int\, \<v,s\>\, d\mu \ \le \ C_{K,v}\|s\|_{W^{1,2}},
    \qquad\text{for all}\quad s\in C^\infty_c(E), \quad \supp s\subset K,
\]
where $\|\cdot\|_{W^{1,2}}$ denotes the norm in the space $W^{1,2}(E)$, cf.
\refe{W12norm}.

\begin{Lem}\label{L:V-u}
For all $u\in W^{1,2}_{\loc}(E)$ we have $V_-u\in W^{-1,2}_{\loc}(E)\cap{}L^1_{\loc}(E)$.
\end{Lem}
\begin{proof}
It follows from \refl{hV-<} by polarization that, for each compact set $K\subset M$,
there exists a constant $C_K>0$ such that
\begin{equation}\label{E:V-ss}
    \Big|\, \int\, \<V_-s_1,s_2\>\, d\mu\, \Big| \ \le \
    C_K\, \|s_1\|_{W^{1,2}}\, \|s_2\|_{W^{1,2}},
    \qquad s_1,s_2\in C^\infty_c(E), \quad\supp s_2\subset K.
\end{equation}

Since the statement of  Lemma is local, we can assume that $u$ is supported in
a coordinate neighborhood. Let us
fix such $u\in W^{1,2}_{\comp}(E)$ and let $u^\rho= \mathcal{J}^{\rho}u$ be as in
\refss{friedrichs}. It follows from \refl{mollifier}, that $u^\rho\to u$ as $\rho\to 0$
both in the space $W^{1,2}(E)$ and in the space $L^2(E)$. In particular,
$\|u^\rho\|_{W^{1,2}}$ is bounded for $0<\rho<1$. Hence, by \refe{V-ss}, for every
compact $K\subset M$ there exists a constant $C_{K,u}>0$, such that
\begin{equation}\label{E:V-u}
 \Big|\, \int\, \<V_-u^\rho,s\>\, d\mu\, \Big| \ \le \
    C_{K,u}\|s\|_{W^{1,2}},
    \qquad s\in C^\infty_c(E), \quad \supp s\subset K,\quad 0<\rho<1.
\end{equation}

Since $u^\rho\to u$ in $L^2_{\loc}(E)$, we have $V_-u^\rho\to V_-u$ in $L^1_{\loc}(E)$.
Hence,
\begin{equation}\label{E:limVurho}
    \lim_{\rho\to 0+}\, \int\, \<V_-u^\rho,s\>\, d\mu
    \ = \
    \int\, \<V_-u,s\>\, d\mu.
\end{equation}
Combining \refe{V-u} and \refe{limVurho}, we obtain
\[
     \Bigl|\int\, \<V_-u,s\>\, d\mu\Bigr| \ \le \  C_{K,u}\|s\|_{W^{1,2}},
    \qquad s\in C^\infty_c(E), \quad \supp s\subset K,
\]
which ends the proof.
\end{proof}

Below we will often denote by $(\cdot,\cdot)$ the duality between $W^{-1,2}_{\comp}(E)$
and $W^{1,2}_{\loc}(E)$. The following lemma shows that this notation will not lead to a
confusion.
\begin{Lem}\label{L:()}
Let $u\in L^2_{\comp}(E)\subset W^{-1,2}_{\comp}(E), \ v\in W^{1,2}_{\loc}(E)$. Then
\[
    (u,v):= \int\, \<u,v\>\, d\mu \ = \ (u,v)',
\]
where  $(\cdot,\cdot)'$ denotes the duality between $W^{-1,2}_{\comp}(E)$ and
$W^{1,2}_{\loc}(E)$ extending the $L^{2}$-scalar product from $\ecomp$ by continuity.
\end{Lem}
\begin{proof}
Using a partition of unity we can assume again that
$u$  and $v$ are supported in a coordinate neighborhood, so we can use the Friedrichs mollifiers.
Let $u^\rho= \mathcal{J}^{\rho}u, v^\rho= \mathcal{J}^{\rho}v$ be as in
\refss{friedrichs}. It follows from \refl{mollifier}, that $u^\rho\to u$ as $\rho\to 0$
both in $W^{-1,2}(E)$ and in  $L^2(E)$. Similarly,
$v^\rho\to v$ in $W^{1,2}(E)$ and in $L^2(E)$. Then
\[
    (u,v) \ = \ \lim_{\rho\to 0+}\, (u^\rho,v^\rho) \ = \
     \lim_{\rho\to 0+}\, (u^\rho,v^\rho)' \ = \ (u,v)'.
\]
\end{proof}

The following lemma is due to H.~Br\'ezis and F.~Browder \cite{Brezis-Browder78}.

\begin{Lem}\label{L:2dualities}
Let $A\in L^2_{\loc}(\End E)$ be a non-negative bundle map
(i.e. for every $x\in M$ the endomorphism $A_x:E_x\to E_x$ has
a non-negative quadratic form),
and let $u\in W^{1,2}_{\comp}(E)$. Assume that $Au\in W^{-1,2}_{\comp}(E)$, so in
fact $Au\in W^{-1,2}_{\comp}(E)\cap L^1_{\comp}(E)$. Then
$\<Au,u\>\in L^1_{\comp}(M)$ and
\begin{equation}\label{E:2dualities}
    \int\, \<Au,u\>\, d\mu \ = \ (Au,u),
\end{equation}
where  $(\cdot,\cdot)$ denotes the duality between $W^{-1,2}_{\comp}(E)$ and
$W^{1,2}_{\loc}(E)$ extending the $L^{2}$-scalar product from $\ecomp$ by continuity.
\end{Lem}

\begin{proof} For the convenience of the reader we  provide a proof, which slightly
differs from the one given in \cite{Brezis-Browder78}.
Using a partition of unity we can assume that $A$ is supported in a coordinate neighborhood.
Since $u$ is only relevant in a neighborhood of $\supp A$, we can assume that
$u$ is supported in the same coordinate neighborhood.
Let $v\in W^{1,2}(E)\cap L^\infty(E)$. We claim that
\begin{equation}\label{E:Auv}
    \int\, \<Au,v\>\, d\mu \ = \ (Au,v).
\end{equation}
Indeed, denote $(Au)^\rho= \mathcal{J}^{\rho}(Au)$. It follows from \refl{mollifier}, that
$(Au)^\rho\to Au$ as $\rho\to 0$ both in  $W^{-1,2}_{\comp}(E)$ and in
$L^1_{\comp}(E)$. Hence,
\[
    (Au,v) \ = \ \lim_{\rho\to 0+}\, \big(\, (Au)^\rho,v\, \big)
    \ = \ \lim_{\rho\to 0+}\,\int\, \big\<\, (Au)^\rho,v\, \big\>\, d\mu.
    \ = \ \int\, \<Au,v\>\, d\mu.
\]

For every $R>0$ define the truncation $u_R$ of $u$ by the formula
\[
    u_R(x) \ = \
    \begin{cases}
    u(x), \qquad&\text{if}\quad |u(x)|\le R;\\
    R\frac{u(x)}{|u(x)|}, \qquad&\text{if}\quad |u(x)| > R.
    \end{cases}
\]
It follows from Theorem~A of the Appendix in \cite{Leinfelder-Simader} that $u_R\in
W^{1,2}_{\comp}(E)$ for all $R>0$ and that $u_R\to u$ as $R\to\infty$ in
$W^{1,2}_{\comp}(E)$. Hence,
\begin{equation}\label{E:AuuR}
    (Au,u) \ = \ \lim_{R\to \infty}\,  (Au,u_R)
    \ = \ \lim_{R\to\infty}\, \int\, \<Au,u_R\>\, d\mu,
\end{equation}
where the last equality follows from \refe{Auv}. By our assumption,
$\<A(x)u(x),u(x)\>\ge 0$ for almost all $x\in M$. Hence,
$\<A(x)u(x),u_R(x)\>$ increases with $R$
for almost all $x\in M$. The lemma follows now from the  monotone
convergence theorem.
\end{proof}

\begin{rem} More general results can be found e.g. in
\cite{Brezis-Browder79}, \cite{Brezis-Browder82}, \cite{Cascante-Ortega-Verbitsky}.
Note that the statement of the lemma is not completely trivial. For example, if
$w, v$ are scalar functions,
$w\in W^{-1,2}_{\loc}\cap L^1_{\loc}$ and $v\in W^{1,2}_{\comp}$,
then it might happen that $wv$ is not in $L^1_{\comp}$,
so that the integral $\int wv\, d\mu$ is not well defined (even though
$(w,v)$ is perfectly well defined from the duality between $W^{-1,2}$ and $W^{1,2}$)
-- see e.g. an example provided in \cite{Brezis-Browder79}.
In fact,  if $w\in L^1_{\loc}$ is fixed, then the condition that $wv$ is in  $L^1_{\loc}$
for every $v\in W^{1,2}_{\comp}$ is equivalent to the inclusion
$|w|\in W^{-1,2}_{\loc}$ (this follows e.g. from \cite{Mazya-book}, Theorem 2 in Sect. 8.4.4).
\end{rem}

\subsection{}
Applying \refl{2dualities} to the bundle map $A=-V_-$ and using \refl{V-u}, we see that
\begin{equation}\label{E:V-uu}
    \int\, \<V_-u,u\>\, d\mu \ = \ (V_-u,u)
\end{equation}
is finite for all $u\in W^{1,2}_{\comp}(E)$. Consider the expression $h_V(u)$ defined by
\begin{equation}\label{EE:formalq} \notag
     h_V(u) \ = \
     \|Du\|^2 \ + \ \int\, \<V_-u,u\>\,d\mu \ + \ \int\, \<V_+u,u\>\, d\mu \ \le \ +\infty,
     \qquad u\in W^{1,2}_{\comp}(E).
\end{equation}
The main result of this section is the following
\begin{Prop}\label{P:simader}
Assume that $h_V(u)\geq \|u\|^2$ for all $u\in\ecomp$. Then
\begin{equation}\label{E:semquad}
     h_V(u) \ \geq \ \|u\|^2,
\qquad\text{for all}\quad u\in\wotcomp(E).
\end{equation}
\end{Prop}
The proof of the proposition occupies the rest of this section.

\medskip
The following simple ``integration by parts" lemma follows e.g. from Theorem 7.7 in
\cite{sh5}.

\begin{lem}\label{L:Lesch}
The equality
\[
    (Du,v) \ = \ (u, D^*v)
\]
holds if one of the sections $u, v$ has compact support and $u\in L^{2}_{\loc}(E)$, $v\in
W^{1,2}_{\loc}(F)$ or, vice versa,  $u\in W^{1,2}_{\loc}(E)$, \, $v\in L^{2}_{\loc}(F)$.
Here $(\cdot,\cdot)$ is understood as either the scalar product in $L^2$ or the duality
between $W^{-1,2}_{\comp}$ and $W^{1,2}_{\loc}$ extending the $L^{2}$-scalar product from
$\ecomp$ by continuity.
\end{lem}

The following well-known lemma (cf. M.~Gaffney~\cite{ga}), whose proof we reproduce for
completeness, provides us with a sufficient amount of ``cut-off" functions to be used
later.
\begin{lem}\label{L:finsler}
Suppose that $g$ is a complete Riemannian metric on a manifold $M$. Then there exists a
sequence of Lipschitz functions $\{\pin\}$ with compact support on $M$ such that
\begin{enumerate}
\item  $0\leq\pin\leq1$ and $|d\pin|\leq\frac{1}{k}$, where $|d\pin|$
denotes the length of the cotangent vector $d\pin$ induced by the metric $g$;
\item $\lim_{k\to\infty}\pin(x)=1$, for any $x\in M$.
\end{enumerate}
\end{lem}
\begin{proof}
Let $d$ be the distance function with respect to the metric $g$. Fix $x_{0}\in M$, and
put $P(x)=d(x,x_{0})$.  Then $P(x)$ is Lipschitz, hence differentiable almost everywhere.
Moreover, $|dP|\leq1$.  The completeness condition $\int^{\infty}ds=\infty$, where $ds$
is the arc-length element associated to $g$, means that $P(x)\to\infty$, as $x\to\infty$.

Consider a function $\chi\in C_{c}^{\infty}(\RR)$ such that $0\leq\chi\leq1$, $\chi=1$
near $0$, and $|\chi'|\leq1$.  Put $\pin(x)=\chi\left(\frac{P(x)}{k}\right)$. Clearly,
$\pin$ has desired properties.
\end{proof}

The following lemma is an analogue of (5.3) from \cite{sh}.
\begin{lem}\label{L:shub}
{}For any $u\in \Dom(H_{V,\max})\subset W^{1,2}_{\loc}(E)$, and any Lipschitz function
with compact support $\psi\colon M\to\RR$
\begin{equation}\label{E:shub}
     h_{V}(\psi u) \ = \ \RE(\psi H_Vu,\psi u)
          \ + \ \|\hat{D}(d\psi)u\|^2.
\end{equation}
\end{lem}
\begin{proof}
Since $u\in W^{1,2}_{\loc}(E)$ we have $D^*Du\in W^{-1,2}_{\loc}(E)$. By \refl{V-u}, we
also have $V_-u\in W^{-1,2}_{\loc}$. Hence, $V_+u= H_Vu- D^*Du- V_-u\in W^{-1,2}_{\loc}$.
Hence, by \refl{2dualities}, we obtain
\[
    \int\, \big\<\, V_+(\psi u),\psi u\, \big\>\, d\mu
    \ = \ \big(\, V_+(\psi u),\psi u\, \big) \ < \ \infty.
\]
In particular,  $h_V(\psi u)<\infty$.
Also, using \refl{()}, we obtain
\begin{multline}\label{E:HVuu=}
    (\psi H_Vu,\psi u) \ = \ \big(\, \psi D^*D u,\psi u\, \big)
    +\big(\, V_+(\psi u), \psi u\, \big)
    +\big(\, V_-(\psi u), \psi u\, \big)
    \\ = \
    \big(\, \psi D^*D u,\psi u\, \big)
        +\int\, \big\<\, V_+(\psi u), \psi u\, \big\>\, d\mu
        + \int\, \big\<\, V_-(\psi u), \psi u\, \big\>\, d\mu.
\end{multline}
Here, in the left hand side $(\cdot,\cdot)$ denotes the $L^2$-scalar product, while
later in this equality it stands for the duality between $W^{-1,2}_{\comp}(E)$ and
$W^{1,2}_{\loc}(E)$.

Using the integration by parts lemma~\ref{L:Lesch}, we get
\begin{multline}\notag
     (D(\psi u),D(\psi u))
    \ = \
     (\hat{D}(d\psi)u,D(\psi u))+ (\psi Du,D(\psi u))
     \ = \
     (\hat{D}(d\psi)u,D(\psi u))+(Du,\psi D(\psi u))
     \\ = \
     (\hat{D}(d\psi)u,\hat{D}(d\psi) u)+
          (\hat{D}(d\psi)u,\psi Du)-(Du,\hat{D}(d\psi)\psi u)
                + (Du,D(\psi^2u))
     \\ = \
     \|\hat{D}(d\psi)u\|^2+ 2i\IM(\hat{D}(d\psi)u,\psi Du)
                         + (\psi D^*Du,\psi u).
\end{multline}
Adding this formula with its complex conjugate and dividing by 2, we obtain
\[
     (D(\psi u),D(\psi u))
     \ = \
     \|\hat{D}(d\psi)u\|^2+ \RE(\psi D^*Du,\psi u).
\]
Now the equality \refe{shub} follows now from \refe{HVuu=}.
\end{proof}
\begin{Prop}\label{P:pospot}
Suppose $H_{V_+}$ (with $V_+\in L^2_{\loc}(\End{}E)$, $V_+\ge 0$) satisfies Assumption B.
(This is true, in particular, if $n\le 3$ or if $D^*D$ has a scalar principal symbol, and
also true if $V_+\in L^p_{\loc}(\End{}E)$ with $p>n/2$ for $n\ge 4$.) Then the operator
$H_{\vp}=D^*D+\vp$ is essentially self-adjoint on $\ecomp$.
\end{Prop}
\begin{proof}
By Theorem X.26 from \cite{rs} it is enough to prove that if $u\in L^{2}(E)$ and
$(H_{\vp}+1)u=0$, then $u=0$.

Indeed, for any such $u$ we have $u\in\Dom(H_{\vp,\max})$. Take $\pik$ as in
\refl{finsler}. By \refl{shub}, with $V=\vp$, we get
\begin{equation}\label{E:equad}
     \hvp(\pik u) \ = \ \RE(\pik H_{\vp}u,\pik u)+\|\hat{D}(d\pik)u\|^2
     \ = \ -\|\pik u\|^2+\|\hat{D}(d\pik)u\|^2.
\end{equation}
Since $\hvp(\pik u)\geq 0$, equation \refe{equad} implies
\begin{equation}\label{E:etquad}
    \|\pik u\|^2 \leq \hvp(\pik u)+\|\pik u\|^2 = \|\hat{D}(d\pik)u\|^2
    \leq\frac{m}{k^2}\|u\|^2,
\end{equation}
where $m=\dim (E_x)$. The last inequality in \refe{etquad} is true by \refl{finsler} and
\refe{equalst}. Taking limit as $k\to\infty$ we get $\|u\|=0$.
\end{proof}

We need the following well known abstract
\begin{lem}\label{L:absform}
Let $\sh$ be a Hilbert space and let $\sd\subset \sh$ be a linear subspace. Let $A:\sd\to
\sh$ be an essentially self-adjoint operator in $\sh$. Suppose $h$ is a closed positive
quadratic form on $\sh$ with domain $Q(h)\supset\sd$ such that
\begin{equation}\label{E:huv}
    h(u,v) \ = \ (Au,v),
    \qquad\text{for all}\quad u\in\sd, \ v\in Q(h),
\end{equation}
where $h(u,v)$ is the sesquilinear form defined from $h$ via polarization. Then the
closure $\oA$ of $A$ is the operator associated to $h$ via the Friedrichs construction,
cf. \cite[\S{}X.23]{rs}. In particular, $\sd$ is the core of $h$.
\end{lem}
\begin{proof}
Let $B$ be the self-adjoint operator associated to $h$ via the Friedrichs construction.
Recall that the domain $\Dom(B)$ consists of the vectors $u\in Q(h)$ such that the map
\begin{equation}\label{EE:lu} \notag
    l_u: v \ \mapsto \ h(u,v), \qquad v\in Q(h),
\end{equation}
is continuous if we endow $Q(h)$ with the norm which it inherits from $\sh$. In this case
$Bu$ is defined from the formula
\begin{equation}\label{E:Bu}
    h(u,v) \ = \ (Bu,v), \qquad v\in Q(h).
\end{equation}

We are ready now to prove the lemma. It follows from \refe{huv} that the map $l_u$ is
continuous for all $u\in \sd$. In other words $\sd\subset \Dom(B)$. Comparing \refe{huv}
and \refe{Bu} we conclude that $Bu= Au$ for all $u\in \sd$. Thus $B$ is a self-adjoint
extension of $A$. Since $A$ is essentially self-adjoint it has a unique self-adjoint
extension and $\oA= B$.
\end{proof}

Let us consider the expression
\begin{equation}\label{EE:posquad} \notag
    h_{V^+}(u) \ = \ \|D u\|^2 \ + \ \int\, \<V_+u,u\>\, d\mu,
\end{equation}
where $h_{V^+}$ is viewed as a quadratic form with the domain
\begin{equation}\nonumber
     Q(h_{V^+}) \ = \ \big\{\, u\in L^2(E)\cap W^{1,2}_{\loc}(E):\,
         h_{V^+}(u)<+\infty\, \big\}.
\end{equation}
Clearly, $h_{V^+}$ is a closed positive form. Recall from \refp{pospot} that the operator
$H_{V^+}$ is essentially self-adjoint on $C^\infty_c(E)$.

\begin{Prop}\label{P:posquadform}
The closure of the operator $H_{V^+}$ is the operator associated to the form $h_{V^+}$
via the Friedrichs construction. In particular, the space $C^\infty_c(E)$ is the core of
the form $h_{V^+}$, i.e. $Q(h_{V^+})$ is the closure of $C^\infty_c(E)$ with respect
to this norm
\[
     \|u\|_1^2 \ : = \ \|Du\|^2+\int\, \<V_+u,u\>\, d\mu +\|u\|^2\,.
\]
\end{Prop}
\begin{proof}
Since $Q(h_{V^+})\subset W^{1,2}_{\loc}(E)$, it follows from the integration by parts
lemma~\ref{L:Lesch}, that
\[
    (H_{V^+}u,v) \ = \ (D^*Du,v)+ \int\, \<V_+u,v\>\, d\mu
    \ = \ (Du,Dv)+\int\, \<V_+u,v\>\, d\mu \ = \ h_{V^+}(u,v),
\]
for all $u\in C^\infty_c(E)$, $v\in Q(h_{V^+})$. Hence, applying Lemma \ref{L:absform}
ends the proof.
\end{proof}

{}For any compact set $K\subset M$ define
\[
    V_-^{(K)}(x) \ = \
    \begin{cases}
         \vm(x), \quad &x\in K\\
         0, \quad &\text{otherwise}
     \end{cases},
     \qquad\qquad\quad V^{(K)} \ = \ V_-^{(K)}+\vp.
\]
Note that $\left|\int\<V_-^{(K)}u,u\>d\mu\right|= |(V_-^{(K)}u,u)|<\infty$ for all $u\in
W^{1,2}_{\loc}(E)$ by \refe{V-uu}.

\begin{lem}\label{L:conver}
Suppose $v_k\to v$ in $\|\cdot\|_1$. Then
 \(
    \lim_{k\to\infty}(V_-^{(K)}v_k,v_k) \ = \ (V_-^{(K)}v,v).
 \)
\end{lem}
\begin{proof}
By definition of $\|\cdot\|_1$, we have
\begin{equation}\label{E:conver2}
    \lim_{k\to\infty}\|D(v-v_k)\| \ = \ 0, \qquad
    \lim_{k\to\infty}\|v-v_k\| \ = \ 0,
\end{equation}
where $\|\cdot\|$ means the $L^2$ norm. Let us fix $\phi\in C^\infty_c(M)$, such that
$\phi=1$ on a neighborhood of $K$. Then, from \refe{conver2} we obtain
\[
    \lim_{k\to\infty}\, \|D(\phi v- \phi v_k)\|
    \ \le \ \lim_{k\to\infty}\, \|\hatD(d\phi)(v-v_k)\| +
    \lim_{k\to\infty}\, \|\phi D(v-v_k)\| \ = \ 0.
\]
In other words, $\phi{v_k}$ converges to $\phi{v}$ in $W^{1,2}_{\comp}(E)$. It follows
now from \refl{hV-<} (or, more precisely, from the equation \refe{V-ss}) that
\[
    \lim_{k\to\infty}(V_-^{(K)}v_k,v_k) \ = \
    \lim_{k\to\infty}(V_-^{(K)}\phi v_k,\phi v_k) \ = \
     (V_-^{(K)}\phi v,\phi v) \ = \ (V_-^{(K)} v, v).
\]
\end{proof}

\subsection{Proof of \refp{simader}}
{}Fix  $u\in W^{1,2}_{\comp}(E)$ and let $K\subset M$  be a compact set, which contains a
neighborhood of $\supp u$. Define
\[
     h_{V^{(K)}}(u) \ = \
     \|Du\|^2 \ + \ \int\, \<\vp u,u\>\, d\mu
        \ + \ \int\, \<V_-^{(K)} u,u\>\, d\mu  \ = \
     \hvp(u) \ + \ (V_-^{(K)}u, u).
\]
If $h_V(u)=+\infty$ then the inequality \refe{semquad} is tautologically true. Hence, we
can assume that $h_V(u)<+\infty$. In particular, $u\in Q(\hvp)$.  By \refp{posquadform},
there exists a sequence $u_k\in\ecomp$ such that $\|u_k-u\|_1\to 0$ as $k\to\infty$.
Using \refe{HV>} and \refl{conver},
we obtain
\begin{multline}\tag*{$\square$}
     \hvp(u)+(V_-^{(K)} u, u)
    \ = \ \lim_{k\to\infty}\hvp(u_k)+ \lim_{k\to\infty}(V_-^{(K)} u_k,u_k)
    \\ = \ \lim_{k\to\infty}h_{V^{(K)}}(u_k) \
    \ \ge \ \lim_{k\to\infty}h_V(u_k)
    \ \ge \ \lim_{k\to\infty}\|u_k\|^2 \ = \ \|u\|^2.
\end{multline}

\begin{cor}\label{C:needinmain}
Let $\delta$, $q$ and $V$ be as in \reft{main}. Then
\[
    \delta(Du,Du)+(Vu,u) \ \geq \ -(qu,u)
    \qquad \text{for all}\quad u\in\wotcomp(E).
\]
\end{cor}
\begin{proof}
It is enough to prove that
\begin{equation}\nonumber
     (Du,Du)+\delta^{-1}\big(\, (V+q)u,u\, \big) \ \geq \  0,
     \qquad\text{for all}\quad u\in W^{1,2}_{\comp}(E).
\end{equation}
This follows immediately from \refp{simader}
\end{proof}

\section{Proof of \reft{main}}\label{S:proof}

Throughout this section we assume that $V$ satisfies the conditions of \reft{main}.

\begin{lem}\label{L:axil}
Let $0\leq\psi\leq q^{-1/2}\le 1$ be a Lipschitz function with compact support and set
\[
     C \ = \sqrt{m}\; \underset{x\in M}\esssup |d\psi(x)|
         \ := \ \underset{x\in M}\esssup
     \left(\RE\, \Tr\left(\left(\symd(d\psi)\right)^{*}\symd(d\psi)\right)\right)^{1/2}.
\]
Then, for any $u\in\Dom(H_{V,\max})$ we have
\begin{equation}\label{E:toprove}
    \|\psi Du\|^{2}
    \ \leq \
    \frac{2}{1-\delta}
     \left(\left(1+\frac{2C^2(1+\delta)^2}{1-\delta}\right)\|u\|^{2}
                                    +\|u\|\|H_Vu\|\right),
\end{equation}
where $\|\cdot\|$ denotes the $L^2$-norm.
\end{lem}
\begin{proof}
{}First note that by \refr{equivalence},
\[
     \underset{x\in M}\esssup|\symd(d\psi)| \ \leq \ C.
\]
Let $u\in \Dom(H_{V,\max})$. Then by Assumption~B, $u\in\wl$. Hence, $\psi^2Du\in
L^2_{\comp}(F)$.  This allows application of the integration by parts \refl{Lesch}, which
implies
\begin{multline}\notag
     \|\psi Du\|^2 \ = \ (D^*(\psi^2Du),u)
     \ = \
     (\psi^{2}D^*Du,u) + 2(\psi\widehat{D^*}(d\psi)Du,u)
     \\ = \
     \RE(\psi^2D^*Du,u) \ + \ 2\RE(\psi\widehat{D^*}(d\psi)Du,u)
     \ \le \ \RE(\psi^2D^*Du,u) \ + \ 2C\|\psi Du\|\|u\|.
\end{multline}
{}Furthermore,
\begin{multline}\notag
    \RE(\psi^2D^*Du,u)
     \ \ge \
     (\psi Du,\psi Du) \ - \ 2C\|\psi Du\|\|u\|
     \\ = \
     (D(\psi u),D(\psi u)) \ - \ (\hat{D}(d\psi)u,\psi Du)
         \ - \ (\psi Du,\hat{D}(d\psi)u) \ - \ 2C\|\psi Du\|\|u\|
    \\ \ge \
     (D(\psi u),D(\psi u)) \ - \ 4C\|\psi Du\|\|u\|.
\end{multline}
Hence, by~\refc{needinmain},
\begin{multline}\label{EE:1-del} \notag
     (1-\delta) \|\psi Du\|^2
     \ \le \ (1-\delta)\RE(\psi^2D^*Du,u) \ + \ 2C(1-\delta)\|\psi Du\|\|u\|
     \\ = \
     \RE(\psi^2H_Vu,u) \ - \ \delta\RE(\psi^2D^*Du,u) \ - \ (\psi^2Vu,u)
                 \ + \ 2C(1-\delta)\|\psi Du\|\|u\|
     \\ \le \
     \|H_Vu\|\|u\| \ - \  \delta(D(\psi u),D(\psi u))-(V(\psi u),\psi u)
                 \ + \ 2C(1+\delta)\|\psi Du\|\|u\|
     \\ \le \
     \|H_Vu\|\|u\| \ + \ (q\psi u,\psi u)
            \ + \ 2C(1+\delta)\|\psi Du\|\|u\|
    \\ \le \
     \|H_Vu\|\|u\| \ + \ \|u\|^2 \ + \ 2C(1+\delta)\|\psi Du\|\|u\|.
\end{multline}
Using the inequality $2ab\leq\varepsilon a^{2}+\frac{b^{2}}{\varepsilon}$ $(a,b\in\RR)$,
with $\varepsilon =\frac{1-\delta}{2(1+\delta)^2}$, we obtain
\[
     2(1+\delta)C\|\psi Du\|\|u\|
     \ \leq \ \frac{1-\delta}2\|\psi Du\|^{2}+
        2C^{2}\frac{(1+\delta)^2}{1-\delta}\|u\|^{2}.
\]
Therefore,
\[
    (1-\delta)\|\psi Du\|^{2}
     \ \leq \
     \frac{1-\delta}2\|\psi Du\|^{2}
       +\Big(\, 1+\frac{2C^2(1+\delta)^2}{1-\delta}\, \Big)\|u\|^{2}
                +\|u\|\|H_Vu\|,
\]
and~\refe{toprove} immediately follows.
\end{proof}

\begin{lem} \label{L:prediduschaya}
Suppose that the metric $\gtm$ is complete (as described in \refss{Fmetric}). Let $q$ be
as in \reft{main}, and let $u\in \Dom(H_{V,\max})$. Then $q^{-1/2}Du \in \fl$ and
\begin{equation}\label{E:ineqt}
    \|q^{-1/2}Du\| \ \leq \
    \frac{2}{1-\delta}
      \left(\left(1+\frac{2L^2(1+\delta)^2}{1-\delta}\right)\|u\|^{2}
                    +\|u\|\|H_Vu\|\right).
\end{equation}
\end{lem}

\begin{proof}
Using a sequence $\pin$ of Lipschitz functions from~\refl{finsler}, define
$\psi_{k}=\pin\cdot q^{-1/2}$.  Then $0\leq\psi_{k}\leq q^{-1/2}$, and
$|d\psi_{k}|\le|d\pin|\cdot q^{-1/2}+\pin|dq^{-1/2}|$. Therefore,
$|d\psi_{k}|\leq\frac{1}{k}+L$, where $L$ is the Lipschitz constant of $q^{-1/2}$.  Since
$\psi_{k}(x)\to q^{-1/2}(x)$ as $k\to\infty$, the dominated convergence theorem applied
to \refe{toprove} with $\psi=\psi_k$ immediately implies~\refe{ineqt}.
\end{proof}

\subsection{Proof of \reft{main}}\label{SS:prmain}
Let $u,v\in \Dom(H_{V,\max})$, and let $\phi\geq0$ be a Lipschitz function with compact
support. By \refr{rmk} the metric $g^{TM}$ is complete. From \refl{prediduschaya}, we see
that $q^{-1/2}Du$ and $q^{-1/2}Dv$ belong to $\fl$.

By Assumption~B, the sections $u$ and $v$ belong to $\wl$. Hence, $Du, Dv\in
L^2_{\loc}(F)$. Using \refl{Lesch}, we obtain
\begin{eqnarray}
    (\phi u, D^*Dv) \ = \ (D(\phi u),Dv)
         \ &=& \
    (\symd(d\phi)u,Dv)+(\phi Du,Dv)\notag\\
     (D^*Du,\phi v) \ = \ (Du,D(\phi v))
         \ &=& \ (Du,\symd(d\phi)v)+(\phi Du,Dv)\notag
\end{eqnarray}
By \refe{equalst},
  \(\displaystyle
     \underset{x\in M}\esssup|\symd(d\phi)|\leq\sqrt{m}\;\displaystyle
\underset{x\in M}\esssup|d\phi(x)|,
  \)
where $|d\phi|$ denotes the length of the covector $d\phi$ in $g^{TM}$. Therefore,
\begin{multline}\label{E:pravlev}
      |(\phi u,H_Vv)-(H_Vu,\phi v)|
     \ \leq \ |(\symd(d\phi)u,Dv)|+|(Du,\symd(d\phi)v)|
     \\ \leq \ \sqrt{m}\; \underset{x\in M}\esssup \left(|d\phi|q^{1/2}\right)\cdot
            \left(\|u\|\|q^{-1/2}Dv\|+\|v\|\|q^{-1/2}Du\|\right).
\end{multline}

Consider a metric $g:=q^{-1}g^{TM}$. Denote the length a covector in this metric by
$|\cdot|_{g}$. By condition (iii) of the theorem this metric is complete.  Using this
metric, take a sequence of Lipschitz functions $\{\pin\}$ as in \refl{finsler}. Since
$|d\pin|_{g}=q^{1/2}|d\pin|$, we obtain $q^{1/2}|d\pin|\leq\frac{1}{k}$.  Therefore,
$\underset{x\in M}\esssup(|d\pin|q^{1/2}(x))\leq\frac{1}{k}$. Using \refe{pravlev}, we obtain
\[
     |(\pin u,H_Vv)-(H_Vu,\pin v)|
     \ \leq \ \frac{\sqrt{m}}{k}\left(\|u\|\|q^{-1/2}Dv\|+\|q^{-1/2}Du\|\|v\|\right)
        \ \to \  0, \quad\text{as} \quad k\to\infty.
\]
On the other side, by the dominated convergence theorem we have
\[
     (\pin u,H_Vv)-(H_Vu,\pin v) \ \to \ (u,H_Vv)-(H_Vu,v)
     \qquad \text{ as }\quad k\to\infty.
\]
Thus $(H_Vu,v)=(u,H_Vv)$, for all $u, v\in \Dom(H_{V,\max})$. Therefore, $H_V$ is
essentially self-adjoint.\hfill$\square$

\section{Proof of \reft{bound}}\label{S:mpf}

We proceed as in \S5.1 of \cite{sh}. Since $H_V$ is semi-bounded below, there exists a
constant $C>0$, such that $H_V\ge -C$ on $C^\infty_c(E)$.

Adding $(C+1)I$ to $H_V$ we may assume that $H_V\geq I$ on $\ecomp$, i.e.
\begin{equation}\nonumber
    (H_Vu,u)\geq\|u\|^2, \qquad u\in\ecomp.
\end{equation}

It is well known (see, e.g., Theorem~X.26 from \cite{rs}) that $H_V$ is essentially
self-adjoint if and only if the equation $H_Vu=0$ has no non-trivial solutions in
$L^2(E)$ (understood in the sense of distributions). Assume that $H_Vu=0, \ u\in L^2(E)$.
In particular, $u\in\Dom(H_{V,\max})$.  By hypotheses, we know that $u\in
W^{1,2}_{\loc}(E)$.

Let $\phi_k$ be the sequence of Lipschitz compactly supported functions from
\refl{finsler}. Then, $\phi_ku$ belongs to $W^{1,2}_{\comp}(E)$, for any $k=1,2,\ldots$.

Using \refp{simader}, \refl{shub}, and the equality $H_Vu=0$, we obtain
\[
     \|\pik u\|^2 \ \le \ h_{V}(\pik u)
     \ = \ \|\hat{D}(d\pik) u\|^2 \ \le \ \frac{m}{k^2}\|u\|^2.
\]
The last inequality is true by \refl{finsler} and \refe{equalst}. Since,
$\|u\|=\lim_{k\to\infty}\|\pik u\|$, we conclude that $u=0$. \hfill$\square$




\appendix
\section{Friedrichs mollifiers}\label{S:frmoll}

Let $u\in\lolocv$ be a vector valued function on $\RR^n$. Let $\sj$ be an integral
operator whose integral kernel is $\kjr(x-y)\Id$, where $\Id$ is $m\times m$ identity
matrix and $\kjr(x-y)$ is as in \refss{friedrichs}. We define $\ur=\sj u$.

The main result of this appendix is the following proposition, which generalizes results
of K.~Friedrichs \cite{fri} and T.~Kato \cite[Sect.~5, Lemma~2]{Kato72}  to our vector valued
setting.

\begin{Prop}\label{P:socommutator}
Consider a  second-order differential operator
\begin{equation}\label{E:L=}
     L \ = \  \sum_{i,k}\pari\,a_{ik}(x)\park+\sum_{i}b_i(x)\pari+c(x),
\end{equation}
where $a_{ik}(x)$, $b_i(x)$ and $c$ are $m\times m$-matrices, such that the matrix elements
of $a_{ik}(x)$ and $b_i(x)$ are locally Lipschitz functions in $\RR^n$, and
the matrix elements of
$c(x)$ belong to $L^\infty_{\loc}(\RR^n)$. Let us suppose that $u\in\lolocv$ and
$Lu\in\lolocv$. Assume, in addition, that
\begin{equation}\label{E:regular}
    \park u \in \lolocv \quad \text{if there exist $x\in \RR^n, \ i\in \{1\nek n\}$ such
    that $a_{ik}(x)\not=0$.}
\end{equation}
Then $L\ur\to Lu$ in $\lolocv$.
\end{Prop}
\begin{rem}\label{R:regular}
Note that the assumption \refe{regular} is tautologically true if $L$ is a first order
operator, i.e., if $a_{ik}\equiv 0$ for all $j,k$. On the other extreme, if $L$ is an
elliptic second order differential operator with smooth coefficients, then \refe{regular}
is a consequence of our other assumptions $u\in\lolocv$ and $Lu\in\lolocv$. This follows
from standard elliptic regularity results, cf., e.g., the arguments at the end of the
proof of \refl{VuL2}. In fact, in this case it is enough to assume that $a_{ik}(x)$ and
$b_i(x)$ are locally Lipschitz  and $c(x)$ is in $L^\infty_{\loc}$ (cf. sketch of the proof
in \cite[Sect.~5,~Lemma~2]{Kato72}\footnote{It is assumed in \cite{Kato72} that
$a_{ik}$ and $b_i$ are in $C^1$ but the same argument can be carried through if we only
assume that they are Lipschitz.}).
\end{rem}

In the proof we will use the following  version of a Friedrichs result (cf. equation (3.8)
in the proof of Main Theorem in \cite{fri}).

\begin{Prop}\label{P:friedrichs}
Let $\jr$ be as in \refe{avg}. Fix $i\in \{1\nek{}n\}$ and let $T=b(x)\pari$, where
$b(x)$ is a locally Lipschitz function in $\RR^n$. Assume that $v\in\loloc(\RR^n)$. Then, as
$\rho\to 0+$, we have
\begin{equation}\label{E:friedrichs}
    (\jr T-T\jr)v\to 0 \qquad \text{in}\quad L^1_{\loc}(\RR^n).
\end{equation}
The same holds if we replace $T=b(x)\pari$ by $\tilde T=\pari\cdot b$, where $b$ is understood
as the multiplication operator by $b(x)$.
\end{Prop}
%
\subsection{Proof of \refp{friedrichs}}\label{SS:prfriedrichs}
Since the statement of \refp{friedrichs} is local, it is enough to prove it in the case
when the support of $b$ is compact and $v\in L^1(\RR^n)$, which we will henceforth
assume. Then $K^\rho:= J^\rho T-TJ^\rho$
 is a continuous operator $L^1(\RR^n)\to
L^1(\RR^n)$. Let $\|K^\rho\|$ denote its norm.

We will use the fact that any locally  Lipschitz function is differentiable almost
everywhere and its pointwise derivative coincides with its distributional derivative
(in particular, the derivative is in  $L^\infty_{\loc}$). This allows to
do integration by parts with Lipschitz functions in the same way as if they
are in $C^1$ (see e.g. \cite{Mazya-book}, Sect. 1.1 and 6.2).

We will start with establishing an analogue of Proposition \ref{P:friedrichs}
for zero order operators.

\begin{lem}\label{L:zero}
If \/ $c\in L^\infty_{\loc}(\RR^n)$,
and $u\in \loloc(\RR^n)$, then \/
$cu^\rho- (cu)^\rho\to 0$ in $\loloc(\RR^n)$ as $\rho\to 0+$.
The same holds  if $u$ is a vector-function (with values in $\CC^m$) and $c$
has values in $m\times m$-matrices.
\end{lem}

\begin{proof}
By \refl{mollifier}(ii), both summands in the right hand side of the equation
\[
    cu^\rho-(cu)^\rho \ = \ c\, \big(\, \sj u-u\, \big) \ + \ \big(\, cu-\sj(cu)\, \big)
\]
converge to 0 in $\lolocv$ as $\rho\to 0+$.
\end{proof}

\begin{lem}\label{L:kernel}
$K^\rho:L^1(\RR^n)\to L^1(\RR^n)$ is the integral operator with Schwartz kernel
\begin{equation}\label{E:kernel}
     \kkr(x,y) \ = \ \deryv\big(\, (b(x)-b(y))\kjr(x-y)\, \big).
\end{equation}
\end{lem}
\begin{proof}
{}For all $v\in C^\infty_c(\RR^n)$ we have
\begin{multline}\nonumber
    (T\jr v)(x)
    \ = \ b(x)\, \derxv\int\, \kjr(x-y)\, v(y)\,dy
    \ = \ \int\,  b(x)\, \derxv\big(\, \kjr(x-y)\, \big)\,v(y)\,dy
    \\ \ = \
    \int\, b(x)\, \left(-\deryv\kjr(x-y)\right)v(y)\,dy
    \ = \ \int \left(-\deryv\big(\, b(x)\kjr(x-y)\,\big)\right)\,v(y)\,dy.
\end{multline}
\begin{equation}\nonumber
    (\jr T v)(x)
    \ = \ \int\, \kjr(x-y)\, b(y)\, \deryv v(y) \, dy
    \ = \ \int \left(-\deryv\big(\, b(y)\, \kjr(x-y)\,\big)\right)\, v(y)\, dy.
\end{equation}
Hence,
\begin{equation}\label{EE:Krho=krho} \notag
    K^\rho v(x) \ = \ \int \, k_\rho(x,y)\, v(y)\, dy,
\end{equation}
where $k_\rho$ is given by \eqref{E:kernel}.
Since $K^\rho$ is a continuous operator $L^1(\RR^n)\to L^1(\RR^n)$ the lemma is proven.
\end{proof}

\begin{lem}\label{L:Krho<C}
There exists a constant $C>0$, such that
\begin{equation}\label{E:k<C}
    \int\, |k_\rho(x,y)|\, dx \ \le \ C, \qquad\text{for all}\quad y\in \RR^n, \, \rho>0.
\end{equation}
In particular, $\|K^\rho\|\le C$ for all $\rho>0$.
\end{lem}
\begin{proof}
{}From \refe{kernel}, we obtain
\begin{equation}\label{E:splitker}
     k_\rho(x,y)
     \ = \
     -\frac{\partial b(y)}{\partial y^i}\,\kjr(x-y)+\big(\, b(x)-b(y)\, \big)\, \deryv\kjr(x-y).
\end{equation}
We now estimate the integral of the absolute values of all terms in the
right-hand side of \refe{splitker}. Since $\int \kjr(x-y)\,dx=1$, we obtain
\begin{equation}\label{E:fst0}
     \int\, \Big|\, \frac{\partial b(y)}{\partial y^i}\,\kjr(x-y)\,\Big|\,dx
     \ \leq \ \underset{y\in\supp b}{\rm ess\,sup}\, \Big|\, \frac{\partial
b(y)}{\partial y^i}\,\Big|
     \ \leq \   \underset{y\in\supp b}{\rm ess\,sup}\, \big|\, \n b(y)\,\big|.
\end{equation}
Recall that $j_\rho(x)=\rho^{-n}j(\rho^{-1}x)$, cf. \refss{friedrichs}.
Hence,
we obtain
\begin{equation}\label{E:scd}
     \int\, \Big|\, \big(\, b(x)-b(y)\, \big)\, \deryv\kjr(x-y)\, \Big|\, dx
     \ \leq \
     \Big(\,  \rho\,\underset{\xi\in\supp b}{\rm ess\,sup}\,
       |\n b(\xi)|\, \Big)
     \cdot\rho^{-1}\, \int \Big|\frac{\partial j(x)}{\partial x^i}\Big|\, dx.
\end{equation}
{}From \refe{fst0} and \refe{scd} we see that \refe{k<C} holds with \/
 \(\displaystyle
    C  =  \Big(\,1 + \int |\n j(x)|\, dx\, \Big)
    \underset{y\in\supp b}{\rm ess\,sup}\, \big|\, \n b(y)\,\big|.
 \)
\end{proof}

Recall that we assume that the support of $b$ is compact. The following lemma summarizes
some additional properties of $k_\rho(x,y)$.
\begin{lem}\label{L:krho}
\begin{enumerate}
\item Support of $k_\rho$ is contained in the $\rho$-neighborhood of\/
$\supp b\times\supp b\subset {\RR^n\times\RR^n}$;\label{1}
\item $k_\rho(x,y)=0$ if $|x-y|>\rho$;\label{3}
\item there exists a constant $C_1>0$, such that $\int |k_\rho(x,y)|dx\,dy\le C_1$ for all
$\rho>0$;\label{2}
\item  $\int{}k_\rho(x,y)\,dy=0$ for all $x\in\RR^n$.\label{4}
\end{enumerate}
\end{lem}
\begin{proof}
(\ref{1}) and (\ref{3}) follow immediately from \refe{kernel}. (\ref{2}) follows from
(\ref{1}) and \refe{k<C}. (\ref{4}) holds because $k_\rho$ is the derivative of the
compactly supported function $\big(b(x)-b(y)\big)j_\rho(x-y)$.
\end{proof}

We now return to the proof of \refp{friedrichs}. In view of \refl{Krho<C} it is enough to
prove \refe{friedrichs} for $v\in C^\infty_c(\RR^n)$. Using Lemmas~\ref{L:Krho<C} and
\ref{L:krho} we obtain
\begin{multline}\notag
    \int\, \big|\, (\kr v)(x)\,\big|\, dx
    \\ = \
    \int \, \Big|\, \int\, \kkr(x,y)\, v(y)\,dy\, \Big|\, dx
    \ = \
    \int\, \Big|\, \int\, \kkr(x,y)\, \big(\, v(y)-v(x)\, \big)\,dy\, \Big|\,dx
    \\ \leq \ C_1\, \max_{|x-y|\le\rho}\, |v(y)-v(x)| \  \le \
    C_1\rho\max_{\xi\in \supp v}\,
      \big|\,\n v(\xi)\, \big| \ \to \ 0, \qquad\text{as}\quad \rho\to 0+.
\end{multline}
This proves the first part of Proposition \ref{P:friedrichs}. To prove the second part
note that
$$
\tilde T u=\pari(b u)=T u + (\pari b) u,
$$
so the desired result follows from the first part and Lemma \ref{L:zero}.
\endproof

\subsection{Proof of \refp{socommutator}}
By (ii) in \refl{mollifier}, $(Lu)^{\rho}-Lu\to 0$ in $\lolocv$ as $\rho\to 0+$. Thus it
is enough to show that
\begin{equation}\label{E:sumrule}
    L\ur-(Lu)^{\rho} \ = \ \big(\, L\ur-Lu\, \big) \ - \ \big(\, (Lu)^\rho-Lu\, \big)
    \ \to 0, \quad \text{in} \ \lolocv \ \text{as} \ \rho\to 0+.
\end{equation}
We will estimate the left hand side of \refe{sumrule} separately for zero, first, and
second order terms of the operator $L$, cf. \refe{L=}. For the zero order terms
we can use Lemma \ref{L:zero}. Proposition \ref{P:friedrichs} gives us desired result
for the first order terms (for vector functions we should apply it separately
to each component of the vector). The following Lemma takes care of the second order terms.

\begin{lem}\label{L:second}
{}Fix \/ $i,k\in \{1\nek n\}$ and \/ $u\in \lolocv$ such that \/ $\park{}u\in \lolocv$.
Let $a(x)$ be an \/ $m\times{}m$-matrix with locally Lipschitz coefficients. Then \/
$\pari(a\park{}u^\rho)-(\pari(a\park u))^\rho\to 0$  in $\lolocv$ as $\rho\to 0+$.
\end{lem}
\begin{proof}
Since $\sj$ commutes with $\park$, we have
\begin{equation}\label{E:commutator}
    \pari(a\park\sj u) - \sj (\pari(a\park u)) \ = \
       \pari(a\sj \park u) - \sj(\pari(a\park u)) \ = \
          \pari(a\sj v)-\sj(\pari(av)),
\end{equation}
where $v=\park u\in\lolocv$.
The lemma follows  now from \refe{commutator} and
Proposition \ref{P:friedrichs}.
\end{proof}

Combining Lemma \ref{L:zero}, Proposition \ref{P:friedrichs} and Lemma \ref{L:second}
we obtain \refp{socommutator}.\hfill$\square$

\subsection{Proof of \refp{keykato}}
Since the statement is local, we work in a coordinate neighborhood $U$, where $\n^*\n $
is given by \refe{katopaper}. Using standard elliptic regularity argument (cf., e.g., the
arguments at the end of the proof of \refl{VuL2}) we see that $u\in W^{1,1}_{\loc}(E)$
i.e. $\pari u\in L^1_{\loc}$ for all $i=1,\dots,n$ in any local coordinates. Hence
\refp{socommutator} immediately gives: $\n^*\n \ur-(\n^*\n u)^{\rho}\to 0$ in
$\loloc(E|_U)$, as $\rho\to 0+$. \hfill$\square$

\section{Positivity}\label{S:positivity}

\def\H{{\mathcal H}}
\def\pa{\partial}
\def\R{\mathbb R}
\def\C{\mathbb C}

The content of this Appendix is mostly well known but we were not able to find an
adequate reference in the literature.

We will briefly describe classical positivity results for the Green function of the
operator $b+\Del_M$ on a complete Riemannian manifold $(M,g)$ with a smooth positive
measure $d\mu$. Here $b>0$ and $\Del_M=d^*d\ge 0$ is the scalar Laplacian, cf.
\refd{Bochner}. Note, that if $d\mu$ is the Riemannian volume form on $M$, then $\Del_M=
-\Del_g$, where $\Del_g$ is the {\em metric} Laplacian, $\Del_gu=\div(\grad{}u)$.

It is a classical fact, due to M.~Gaffney \cite{ga}, that the operator $\Delta_g$ is
essentially self-adjoint. Gaffney's argument  works for $\Del_M$ without any change (for
any measure $d\mu$).  We reproduce it in the proof of  Lemma \ref{L:finsler}. A more
general statement about the essential self-adjointness of $\n^*\n$ follows also from
\refp{pospot}. (Note that Assumption~A is tautologically true in this case, while
Assumption~B follows from elliptic regularity.)

It follows that the operator $b+\Del_M$ is also essentially self-adjoint and strictly
positive. More precisely,
\begin{equation}
((b+\Del_M)u,u)\ge b(u,u), \quad u\in C_c^\infty(M),
\end{equation}
hence this holds for all $u\in\Dom(\Del_{M,\max})$. In particular, the spectrum of
$b+\Del_M$ is a subset of $[b,+\infty)$, hence the operator $(b+\Del_M)^{-1}$ is
everywhere defined and bounded. Denote by $G_b=G_b(\cdot,\cdot)$ its Schwartz kernel,
which is locally integrable  on $M\times M$, smooth ($C^\infty$) off the diagonal in
$M\times M$ and has singularities on the diagonal which are easy to describe by a variety
of methods (e.g. using technique of pseudodifferential operators).

\begin{Thm}\label{T:positivity} In the notations above,
\begin{equation}
G_b(x,y) \ \ge \ 0 \qquad \text{\rm for all} \quad x,y\in M,\ x\ne y.
\end{equation}
\end{Thm}

\begin{proof} This proof was communicated to us by A.~Grigoryan.
It suffices to establish that
\begin{equation}
u(x):=\int_M G_b(x,y)\phi(y)d\mu(y)\ge 0,\quad x\in M,
\end{equation}
for any $\phi\in C_c^\infty(M)$, $\phi\ge 0$. Clearly $u=(b+\Del_M)^{-1}\phi\in L^2(M)$
and it satisfies the  equation
\begin{equation}\label{E:eq-for-u}
(b+\Del_M)\, u \ = \ \phi.
\end{equation}
In particular, $u\in C^\infty(M)$.

Note that the equation \eqref{E:eq-for-u} has a unique solution $u\in L^2(M)$ (understood
in the sense of distributions). Now we will give a construction of a positive solution
$v\in L^2(M)$ of this equation which will end the proof since then we would have $v=u$
due to the uniqueness.

Let us take a sequence of relatively compact open subsets with smooth boundary in $M$
\begin{equation*}
\Omega_1\Subset\Omega_2\Subset\dots\Subset \Omega_k\Subset\Omega_{k+1}\Subset\dots
\end{equation*}
(i.e. $\Omega_k$ is relatively compact in $\Omega_{k+1}$), such that they exhaust $M$,
i.e. their union is $M$. For any $k$ denote by $v_k$ the solution of the following
Dirichlet problem in $\Omega_k$
\begin{equation}\label{E:Dirichlet-k}
  \big(\, b+\Del_M\, \big)\, v_k
  \ = \
  \phi,\qquad v_k|_{\partial \Omega_k}  =  0.
\end{equation}
We will only consider sufficiently large $k$, so that $\supp\phi\subset\Omega_k$. It
follows from the maximum principle that $v_k\ge 0$ and further that
\begin{equation}\label{EE:ineq-k} \notag
0\le v_k\le v_{k+1} \qquad \text{\rm for all}\ \ k.
\end{equation}
So the sequence of functions $v_1,v_2,\dots$ is increasing, hence there exists a
pointwise limit
\begin{equation}\label{EE:v-limit} \notag
  v(x) \ = \ \lim_{k\to\infty}\,v_k(x),\qquad x\in M.
\end{equation}
Let us prove that it is in fact  everywhere finite and locally bounded.

To this end multiply the equation in \eqref{E:Dirichlet-k} by $v_k d\mu$ and integrate
over $\Omega_k$. Using the Stokes formula, we obtain
\begin{equation*}\label{EE:Stokes-k} \notag
   \int_{\Omega_k}\,(bv_k^2+|dv_k|^2)\,d\mu
   \ = \
   \int_{\Omega_k}\, \phi v_k\, d\mu
   \ \le \
    \frac{b}{2}\,
      \int_{\Omega_k}\, v_k^2\, d\mu \ + \
         \frac{1}{2b}\, \int_{\Omega_k}\,\phi^2\, d\mu,
\end{equation*}
hence
\begin{equation*}\label{EE:L2-vk-est} \notag
   \int_{\Omega_k}\, \left(\,\frac{b}{2}v_k^2+|dv_k|^2\,\right)\,d\mu
   \ \le \
    \frac{1}{2b}\, \int_{\Omega_k}\,\phi^2\, d\mu,
\end{equation*}
which gives, in particular, uniform estimate of the $L^2$ norm of $v_k$ on any compact
set in $M$. Due to the standard interior elliptic estimates
(see e.g. \cite[Sect. 5.3]{Triebel}) this implies that every
derivative of $v_k$ is bounded uniformly in $k$ on every compact set $L\subset M$, hence
the sequence $v_k$ converges in the topology of $C^\infty(M)$. Therefore $v$ is
everywhere finite, positive and satisfies \eqref{E:eq-for-u}. Applying the Fatou theorem
we also see that $v\in L^2(M)$. Therefore $v=u$, hence $u\ge 0$.
\end{proof}

In a more general context, establishing positivity of a solution $u$ for the equation
\eqref{E:eq-for-u}  might present a problem. Namely, assume that
\begin{equation}\label{EE:u-nu} \notag
     \big(\, b+\Del_M\, \big)\, u \ = \ \nu \ \ge \ 0,
     \qquad u\in L^2(M),
\end{equation}
where the inequality $\nu\ge 0$ means that $\nu$ is a positive distribution, i.e.
$(\nu,\phi)\ge 0$ for any $\phi\in C_c^\infty(M)$. It follows that $\nu$ is in fact a
positive Radon measure (see e.g. \cite{Gelfand-Vilenkin}, Theorem 1 in Sect. 2, Ch.II).

\subsection*{Conjecture P} In this situation $u\ge 0$ (almost  everywhere or,
equivalently, as a distribution).

\medskip
At a first glance it seems that the above proof (with $\phi\in C_c^\infty(M)$ instead of
$\nu$) might work in this case as well. But in fact it does not work without additional
restrictions on $u$ or $M$ (see Proposition \ref{P:cut-off} below). Let us clarify what
the difficulty is.

Taking a test function $\phi\in C_c^\infty(M)$, $\phi\ge 0$,  we need to prove that
$(u,\phi)\ge 0$. Let us solve the equation $(b+\Del_M)\psi=\phi$, $\psi\in L^2(M)$. We
see then that $\psi\in C^\infty(M)$ and $\psi\ge 0$ according to Theorem
\ref{T:positivity}. So we can write
\begin{equation*}\label{EE:int-by-parts0} \notag
     (u,\phi) \ = \ \big(\, u,\, (b+\Del_M)\,\psi\, ).
\end{equation*}
Now the right hand side can be rewritten as
\begin{equation}\label{EE:int-by-parts1} \notag
     \big(\, u,\,(b+\Del_M)\,\psi\,\big)
     \ = \
     \big(\, (b+\Del_M)\,u,\psi\, \big)_S \ = \ (\nu,\psi)_S,
\end{equation}
where the hermitian form $(\cdot,\cdot)_S$ in the right hand side is obtained by
extending the scalar product in $L^2(M)$ by continuity (from $C_c^\infty(M)\times
C_c^\infty(M)$) to (non-degenerate) hermitian duality of Sobolev spaces
\begin{equation}\label{E:int-by-parts}
     \widetilde H^{-2}(M)\times \widetilde H^{2}(M) \ \to \ \C,
\end{equation}
where $\widetilde H^{-2}(M)=(b+\Del_M)L^2(M)$, $\widetilde H^2(M)=\Dom
(\Del_{M,\max})=\Dom(\Delta_{M,\min})$ with the Hilbert structures transferred from
$L^2(M)$ by the operator $b+\Del_M$, acting from $L^2(M)$ to $\widetilde H^{-2}(M)$ and
from $\widetilde H^2(M)$ to $L^2(M)$. The norms in $\widetilde H^2(M)$ and $\widetilde
H^{-2}(M)$ are given by the formulas
\begin{equation}\label{E:H2-norm}
      \|v\|_2 \ = \ \big\|\, (b+\Del_M)v\, \big\|,
      \qquad \big\|\, (b+\Del_M)\, f\, \|_{-2} \ = \ \|f\|,
\end{equation}
where $\|\cdot\|$ is the norm in $L^2(M)$.

The extension of the duality \eqref{E:int-by-parts} from $C_c^\infty(M)\times
C_c^\infty(M)$ is well defined because $C_c^\infty(M)$ is dense in both spaces
$\widetilde H^{-2}(M)$ and from $\widetilde H^2(M)$ in the corresponding norms
\eqref{E:H2-norm}. Indeed, density of $C_c^\infty(M)$ in $\widetilde
H^2(M)$ means simply that $(b+\Del_M)C_c^\infty(M)$ is dense in $L^2(M)$. To establish
this, let us take $f\in L^2(M)$ which is orthogonal to $(b+\Del_M)C_c^\infty(M)$ in
$L^2(M)$. This means that $(b+\Del_M)f=0$ in the sense of distributions, i.e. $f$ is in
the null-space of the maximal operator $(b+\Del_M)_{\max}$. This implies that $f=0$ due
to the above mentioned essential self-adjointness and strict positivity of $b+\Del_M$.

Similarly, density of $C_c^\infty(M)$ in $\widetilde H^{-2}(M)$ means that
$(b+\Del_M)^{-1}C_c^\infty(M)$ is dense in $L^2(M)$. To prove this, consider $h\in
L^2(M)$ such that $h$ is orthogonal to $(b+\Del_M)^{-1}C_c^\infty(M)$. Since
$(b+\Del_M)^{-1}$ is a bounded self-adjoint operator, this would imply that
$(b+\Del_M)^{-1}h=0$, hence $h=0$.

Note that the space $\widetilde H^{2}(M)$ is different (at least formally) from the space
$H_2^2(M)$ whose norm includes arbitrary second-order covariant derivatives (see E.~Hebey
\cite{Hebey},  Sect. 2.2). It seems still unknown whether $C_c^\infty(M)$ is dense in
$H_2^2(M)$ (see  Section 3.1 in \cite{Hebey}).

Now we need to know whether it is true or not that
\begin{equation}\label{E:equiv-of-dualities}
(\nu,\psi)_S=\int_M \psi\nu
\end{equation}
(the integral in the right hand side makes sense as the integral of a positive measure,
though it can be infinite.) If this is true then we are done because the integral is
obviously non-negative.

A possible way to establish \eqref{E:equiv-of-dualities} is to present the function
$\psi$ as a limit
\begin{equation}\label{EE:psi-limit} \notag
      \psi \ = \ \lim_{k\to\infty}\, \psi_k,
\end{equation}
where $\psi_k\in C_c^\infty(M)$, $\psi_k\ge 0$, $\psi_k\le\psi_{k+1}$, and the limit is
taken in the norm $\|\cdot\|_2$. The equality
\eqref{E:equiv-of-dualities} obviously
holds if we replace $\psi$ by $\psi_k$, so in the limit we obtain the equality for
$\psi$.

We can try to take $\psi_k=\chi_k\psi$, where $\chi_k\in C_c^\infty(M)$, $0\le\chi_k\le
1$, $\chi_k\le\chi_{k+1}$ and for every compact $L\subset M$ there exists $k$ such that
$\chi_k|_L=1$. Then, obviously, $\psi_k\to\psi$ as $k\to +\infty$ in $L^2(M)$. We also
want to have $\Del_M\psi_k\to\Del_M\psi$ in $L^2(M)$. Clearly
\begin{equation}\label{EE:Delta-product} \notag
     \Del_M \psi_k \ = \ \chi_k\Del_M\psi-2\<d\chi_k,d\psi\> +(\Del_M\chi_k)\psi.
\end{equation}
Obviously, $\chi_k\Del_M\psi\to\Del_M\psi$ in $L^2(M)$. Now note also that
\begin{equation}\label{EE:d-estimate} \notag
     \|d\psi\|^2 \ = \ (\Del_M\psi,\psi)\le
\frac{1}{2}\|\Del_M\psi\|^2+\frac{1}{2}\|\psi\|^2
\end{equation}
and on the other hand $d\chi_k\to 0$ and $\Del_M\chi_k\to 0$ in $C^\infty(M)$. Our goal
will be achieved if we can construct $\chi_k$ in such a way that
\begin{equation}\label{E:chi-k-estimate}
     \sup_{x\in M}|d\chi_k(x)|
     \ \le \ C, \qquad \sup_{x\in M}|\Del_M\chi_k(x)|\le C,
\end{equation}
where $C>0$ does not depend  on $k$. This leads to the following statement, communicated
to us by E.~B.~Davies:

\begin{Prop}\label{P:cut-off}
Let us assume that $(M,g)$ is a complete Riemannian manifold with a positive smooth
measure, such that there exist cut-off functions $\chi_k$, $k=1,2,\dots,$ satisfying the
assumptions above. Then Conjecture P holds.
\end{Prop}

It is always possible (on any complete Riemannian manifold) to construct cut-off
functions $\chi_k$ satisfying all the conditions except the uniform estimate for
$\Del_M\chi_k$ (i.e. the second estimate in
\eqref{E:chi-k-estimate}). This was done by
H.~Karcher \cite{ka} (see also \cite{sh}). But the estimate of $\Del_M\chi_k$ presents a
difficulty. The following sufficient condition provides an important class of examples.

\begin{Prop}\label{P:bounded-geometry}
Let $(M,g)$ be a manifold of bounded geometry with arbitrary measure $d\mu$ having a
positive smooth density. Then Conjecture P holds.
\end{Prop}
\begin{proof} For the definition and properties of manifolds of
bounded geometry see e.g. \cite{Roe, Shubin92}. In particular, a construction of cut-off
functions $\chi_k$ satisfying all the necessary properties on any manifold of bounded
geometry can be found in  \cite{Shubin92}, p.61.
\end{proof}

\begin{rem} Note that the bounded geometry conditions are not
imposed on the measure $d\mu$ which can be an arbitrary measure with a positive smooth
density.
\end{rem}

\subsection{Example}\label{SS:P-on-Rn}
Let us consider $M=\R^n$ with a metric $g$ which coincides with the standard flat metric
$g_0$ outside a compact set. Then Conjecture P holds. It can also be proved by using
invertibility of $b+\Del_M$ in the Schwartz space ${\mathcal S}(\R^n)$, which in turn
implies its invertibility in the dual space ${\mathcal S}'(\R^n)$ consisting of tempered
distributions.

\section{Stummel and Kato classes}\label{S:C}

In this Appendix we will briefly review definitions and the most important properties of
Stummel and Kato classes of functions on $\R^n$ and on manifolds.

\subsection{Stummel classes}\label{SS:Stummel-class}
Uniform  Stummel classes on $\R^n$ were introduced by F.~Stummel \cite{Stummel}. More
details  about Stummel classes and proofs can be found in Sect. 1.2 in \cite{CFKS}, Ch. 5
and 9 in \cite{Schechter}, and also \cite{Aizenman-Simon, Simon00, Stummel}.

The (uniform) Stummel class $S_n$ consists of measurable real-valued functions $V$ on
$\R^n$, such that
\begin{align}\notag
    &\lim_{r\downarrow 0}\,
     \left[\, \sup_x\int_{|x-y|\le r}|x-y|^{4-n}\, |V(y)|^2\, dy\, \right] \ = \ 0
        &\text{\rm if}\quad  n\ge 5;\\ \notag
    &\lim_{r\downarrow 0}\,
      \left[\sup_x\int_{|x-y|\le r}\ln(|x-y|^{-1})\, |V(y)|^2\, dy\, \right] \ = \ 0
         &\text{\rm if}\quad  n=4;\\ \notag
     & \sup_x \int_{|x-y|\le r_0}|V(y)|^2\, dy \ < \ \infty &\text{\rm if}\quad n\le 3.
\end{align}
Here $r_0>0$ is arbitrarily fixed. Clearly, this class $S_n$ is invariant under
multiplication by any real-valued bounded measurable function. It is also easy to see
that if  $f\in S_n$ and a diffeomorphism $\phi:\R^n\to\R^n$ is linear near infinity, then
$\phi^*f=f\circ\phi\in S_n$. Therefore the corresponding local version $S_{n,\loc}(M)$ is
well defined for any manifold with $\dim M=n$. We will denote by $S_{n,\comp}(M)$ the
class of functions $f\in S_{n,\loc}(M)$ which have compact support. Both classes
$S_{n,\loc}(M)$ and $S_{n,\comp}(M)$ are invariant under diffeomorphisms of $M$ and also
under multiplication by any real-valued locally bounded measurable function.

There is the following relation of the local Stummel classes with $L^p$:
\begin{equation}\label{E:Lp-Stummel}
L^p_{\loc}(M)\subset S_{n,\loc}(M), \ \text{if}\ p>n/2\ \text{for} \ n\ge 4; \ p=2\
\text{for} \ n\le 3.
\end{equation}
The same holds for the classes of functions with compact support and for the uniform
classes $S_n$ on $\R^n$ (if $L^p(\R^n)$ is replaced by the class
$L^p_{\text{unif}}(\R^n)$ which consists of functions $u\in L^p_{\loc}(\R^n)$ whose
$L^p$-norms on all unit balls are bounded).

An important property of uniform Stummel classes is that if $u\in S_k$ on $\R^k$ and
$\pi:\R^n\to \R^k$ is a surjective linear projection, then $\pi^*u\in S_n$. This
obviously implies corresponding property for the local Stummel classes with respect to
submersions of manifolds.

The following domination property makes the Stummel classes important in the theory of
Schr\"odinger operators: if $V\in S_n$ on $\R^n$, then for any $a>0$ there exists $C>0$
such that
\begin{equation}\label{E:Stummel-dom}
      \|Vu\| \ \le \ a\|\Del u\|+C\|u\|,\quad u\in C_c^\infty(\R^n),
\end{equation}
where $\Del$ is the flat Laplacian. It follows that if $V\in S_{n,\loc}(M)$, where $\dim
M=n$, then \refe{V<Del} holds for $V$ (substituted for $|V_-|$) due to the fact that
locally all elliptic operators of the same order  have equal strength.

\subsection{Example.}\label{SS:Coulomb-dom}
$1/|x|\in L^2_{\text{unif}}(\R^3)$, hence $1/|x|\in S_3$. It follows that $1/|x_k-x_l|\in
S_{3N}$ in $\R^{3N}$ if $1\le k<l\le N$, $x_k,x_l\in \R^3$ and the points in $\R^{3N}$
are presented in the form $x=(x_1,\dots, x_N)$ with $x_j\in \R^3$. Therefore
\eqref{E:Stummel-dom} implies that the domination relation
\begin{equation}\label{Coulomb-dom}
      \||x_k-x_l|^{-1}u\| \ \le \ a\|\Del u\|+C\|u\|, \quad u\in C_c^\infty(\R^{3N}),
\end{equation}
holds with an arbitrary $a>0$ and $C=C(a)$. Note however that $|x_k-x_l|^{-1}\not\in
L^{3N/2}(\RR^{3N})$ already for $N=2$.

\subsection{Kato classes}\label{SS:Kato-class}
The uniform Kato classes $K_n$ are defined in a way which is  similar to the Stummel
classes $S_n$ and play the same role for domination of quadratic forms as the Stummel
classes play for operator domination. They were introduced by R.~Beals \cite{Beals}  as a
particular case of more general classes from M.~Schechter's paper \cite{Schechter67} (see
also M.~Schechter's book \cite{Schechter} or its first edition which was published in
1971). A year later this class was reinvented by T.~Kato \cite{Kato72} in his famous
paper where the Kato inequality made its debut.

A good introduction to the Kato classes $K_n$ can be found in Section 1.2 of \cite{CFKS}.
They were extensively studied by M.~Aizenman and B.~Simon \cite{Aizenman-Simon}
and B.~Simon \cite{Simon82}.

The (uniform) Kato class $K_n$ consists of measurable real-valued functions $V$ on
$\R^n$, such that
\begin{align}\notag
     &\lim_{r\downarrow 0}\,
      \left[\, \sup_x\int_{|x-y|\le r}|x-y|^{2-n}\, |V(y)|\, dy\, \right] \ = \ 0
         &\text{\rm if}\quad  n\ge 3;\\ \notag
     &\lim_{r\downarrow 0}\,
      \left[\, \sup_x\int_{|x-y|\le r}\ln(|x-y|^{-1})\, |V(y)|\, dy\, \right] \ = \ 0
             &\text{\rm if}\quad  n=2;\\ \notag
     & \sup_x\, \int_{|x-y|\le r_0}|V(y)|\, dy \ < \ \infty &\text{\rm if}\quad n=1.
\end{align}
Here $r_0>0$ is arbitrarily fixed. Again $K_n$ is invariant under multiplication by any
real-valued bounded measurable function, as well as under diffeomorphisms  of $\R^n$
which are  linear near infinity. So the corresponding local class $K_{n,\loc}(M)$ is well
defined for any manifold with $\dim M=n$, and we will denote by $K_{n,\comp}(M)$ the
class of functions $f\in K_{n,\loc}(M)$ which have compact support. Both classes
$K_{n,\loc}(M)$ and $K_{n,\comp}(M)$ are invariant under diffeomorphisms of $M$ and also
under multiplication by any real-valued locally bounded measurable function.

The relation similar to \eqref{E:Lp-Stummel} holds for the Kato classes as well:
\begin{equation}\label{EE:Lp-Kato} \notag
     L^p_{\loc}(M)\subset K_{n,\loc}(M),
     \ \text{if}\ p>n/2\ \text{for} \ n\ge 2; \ p=2\ \text{for} \ n=1.
\end{equation}
The same holds for the classes of functions with compact support and for the uniform
classes $K_n$ on $\R^n$ (if $L^p(\R^n)$ is replaced by the class
$L^p_{\text{unif}}(\R^n)$).

As for the uniform Stummel classes, it is also true that if $u\in K_s$ on $\R^s$ and
$\pi:\R^n\to \R^s$ is a surjective linear projection, then $\pi^*u\in K_n$. This
obviously implies the corresponding property for the local Kato classes with respect to
submersions of manifolds.

The following domination property for quadratic forms  makes the Kato classes important:
if $V\in K_n$ on $\R^n$, then for any $a>0$ there exists $C>0$ such that
\begin{equation}\label{EE:Kato-dom} \notag
      |(Vu,u)| \le \ a(\Del u,u)+C\|u\|^2,\quad u\in C_c^\infty(\R^n).
\end{equation}

The arguments given in Example \ref{SS:Coulomb-dom} can be repeated verbatim with Stummel
classes replaced by Kato classes.

\section{More history}\label{S:rev}

In this Appendix we will describe some related results and provide
bibliographical comments. This section complements
Sect.\ref{SS:hist} where the most recent references were provided.

This review is by no means complete. It would be next to impossible to make it complete.
(For example, MathSciNet, the database of American Mathematical Society, based on
Mathematical Review, lists more than 1000 related papers.) So except several landmark
papers and papers which were explicitly important for our work or closely related to it,
we concentrated on results about  operators on manifolds. A comprehensive review of
self-adjointness results for one-dimensional operators can be found in N.~Dunford and
J.~Schwartz \cite{Dunford-Schwartz}. About the multidimensional case (for operators on
$\RR^n$) the reader may consult M.~Reed and B.~Simon \cite[Ch. X]{rs},
D.~E.~Edmunds and W.~D.~Evans \cite{Edmunds-Evans}, and also review papers
by H.~Kalf, U.-W~Schminke, J.~Walter and R.~W\"ust \cite{Kalf-75}, H.~Kalf \cite{Kalf80}
and B.~Simon \cite{Simon00}.  We tried not to repeat these sources unless this was
relevant for the main text of our paper.

As many good stories in mathematics, this one was initiated by H.~Weyl (1909-1910)
in his pioneering papers on spectral theory of one-dimensional
symmetric singular differential operators (see \cite{Weyl1, Weyl2, Weyl3}, and also
\cite[Ch. IX] {Coddington-Levinson}, \cite[Ch. 13]{Hellwig} or
\cite[Ch. II, Sect. 2]{Levitan-Sargsyan}).
Here ``singular" means that one or two ends of the interval where the operator is considered,
is either infinite,
or has a singularity of the coefficients at this end. Depending on the
asymptotic behavior of solutions of the corresponding ordinary
differential equation at a singular end of the
interval,  where the operator is considered, H.~Weyl classifies the situation
at this end as the case of either {\it limit point} or {\it limit circle}.
In modern terminology (which appeared decades later, after invention of quantum mechanics
and its von Neumann mathematical formulation) the case of limit point on both ends
corresponds to essential self-adjointness.

Among the first authors who wrote about  multi-dimensional
Schr\"odinger operators $H_V=-\Delta+V$ in $\RR^n$ we find
T.~Carleman~\cite{Carleman} and K.~Friedrichs~\cite{Friedrichs}
who independently proved the essential self-adjointness in case
when $V$ is locally bounded and semibounded below.
(Carleman's proof was reproduced in the book by I.M.
Glazman~\cite[Ch.~1, Theorem~34]{Glazman}.)

Moving closer to this day, consider a magnetic Schr\"odinger operator
(on $M=\RR^{n}$) which has the form
\begin{equation}\label{E:oper}
H_{g,b,V}=-\sum_{j,k=1}^{n} (\partial_j-ib_j(x))g^{jk}(x)(\partial_k-ib_k(x))+V(x),
\end{equation}
where the matrix $g^{jk}(x)$ (the metric) is assumed to be real, symmetric, positive definite,
$b_j$ and $V$ are some real valued functions. Additional regularity conditions can be
imposed on $g^{jk}$, $b_j$ and $V$.  We will assume that the coefficients
$g^{jk}$ are sufficiently smooth (at least locally Lipschitz).

We will write $H_{b,V}$ for the operator $H_{g,b,V}$ with the flat metric $g^{jk}$
i.e $g^{jk}=\delta_{jk}$.

Assume first for simplicity that  $b_j\in C^1(\RR^n)$.
We wish the operator $H_{b,V}$ to be naturally defined on $C_c^\infty(\RR^n)$.
Then we need to assume that $V\in L^2_{\loc}(\RR^n)$.

The most dramatic improvement of the result by Carleman and Friedrichs
in case of semi-bounded below potentials $V$ is due to T.~Kato \cite{Kato72}.
Kato established that the requirement of local boundedness of $V$ can be
completely removed, so it is enough to require that $V\in L^2_{\loc}(\RR^n)$
and $V$ is semi-bounded below (even for the magnetic Schr\"odinger operator
$H_{b,V}$ with $b_j\in C^1(\RR^n)$).
What we now call Kato's inequality, first appeared in this context in \cite{Kato72}.
Namely, let us denote by $L$ the operator in \refe{oper} with $V=0$, and by $L_0$ the operator
with $b_j=0$ and $V=0$ (so $L_0=-\Delta$). Suppose also that $u\in\lloc^{1}(\RR^n)$
and $Lu\in\lloc^{1}(\RR^n)$.  Then the following inequality holds:
\begin{equation}\label{EE:simpkato} \notag
L_0|u|\leq\RE(Lu\cdot \sign \bar{u}),
\end{equation}
where $(\sign u)(x)=u(x)/|u(x)|$ if $u(x)\neq 0$, and $0$ otherwise.
The main result of the Kato's paper~\cite{Kato72} was the
essential self-adjointness of $H_{b,V}$ with
$V=V_1+V_2$, where $V_1\in\lloc^{2}(\RR^n)$, $V_1\geq-f(|x|)$,  $f(r)$ is a monotone
non-decreasing function of $r=|x|$ which is $o(r^2)$ as $r\to\infty$, and $V_2$
belongs to what is now called Kato's class (cf.  Appendix \ref{S:C}).
In particular, Kato was first to establish essential self-adjointness of the operator
$H_V$ with $V\in L^2_{\loc}(\RR^n)$, $V\ge 0$, without any global conditions on $V$.

Returning to a more general case of the operator $H_{b,V}$,
now omit the requirement  $b_j\in C^1(\RR^n)$.
Then the natural requirements for the minimal operator to be well defined are
$b\in (L^4_{\loc}(\RR^n))^n$, $\div b\in L^2_{\loc}(\RR^n)$ and $V\in L^2_{\loc}(\RR^n)$.
Working under these minimal regularity conditions only, H.~Leinfelder and C.~Simader
\cite{Leinfelder-Simader} were able to improve Kato's result by allowing $V=V_1+V_2$,
where $V_1\ge -c|x|^2$ and $V_2$ is $\Delta$-bounded with a relative bound $a<1$.
In particular, if $V$ is semi-bounded below, then the above minimal
regularity conditions on $b$ and $V$ are sufficient for essential self-adjointness
(see also \cite[Ch. 1]{CFKS}). A non-trivial technique of non-linear truncations
is used in \cite{Leinfelder-Simader} to approximate functions from the maximal
domain of the operator $H_{b,V}$ and its quadratic form
by bounded functions.

A natural question (first formulated probably by I.~M.~Glazman),
is as follows:  can we replace
semi-boundedness below of the potential (in the Carleman-Friedrichs
essential self-adjointness result above)  by semi-boundedness
of the operator itself (i.e. semi-boundedness of the corresponding quadratic form)
on $C_c^\infty(\RR^n)$? (In this case the potential may be not semi-bounded
below on some relatively ``small" sets.)

It occurs that if we only require that $V\in\lloc^{2}(\RR^n)$, then the semiboundedness
below of $H_V$ does not in general imply its essential self-adjointness as we can see
from following example (cf. ~\cite[Ex.~1a]{Kalf-75}, ~\cite[Appendix 2]{Simon} and in
case $n=5$ also \cite[Example 4 in Sect.X.2]{rs}):

\subsection{Example}\label{SS:ex}
Let $H_V=-\Delta+V$, where $V=\frac{\beta}{|x|^2}$ and $n\geq 5$. Then $H_V$ is
essentially self-adjoint on $C_c^{\infty}(\RR^n)$ if and only if
$\beta\geq\beta_0:=1-(\frac{n-2}{2})^2$.  $H_V$ is semibounded below if and only if
$\beta\geq-(\frac{n-2}{2})^2.$ The first assertion follows from separation of variables,
and the second follows from the Hardy inequality (cf. \cite{Kalf-Walter}).

The first result about essential self-adjointness of
semi-bounded differential operators is probably due to F.~Rellich \cite{Rellich51}
who considered one-dimensional case and general Sturm-Liouville operator
on an interval with one or two singular ends.
In case when $M=\RR^n$, $V$ continuous, and  $H_V=-\Delta+V$ is semibounded below, the
essential self-adjointness of $H_V$ was  proven by
A.~Ya.~Povzner~\cite[Th.~6, Ch.~I]{Povzner}. (Povzner writes that the result was conjectured by
I.~M.~Glazman in a conversation between them. This happened not later than 1952 when the Povzner
paper was submitted. It is not clear whether Glazman knew about the existence of Rellich's paper.)
Approximately 1.5 years  after Povzner's paper  appeared in print,  the question about
the self-adjointness of semi-bounded Schr\"odinger  operator
$H_V$ was asked by F.~Rellich \cite{Rellich56} in his Amsterdam talk  at the International Congress of
Mathematicians.  Answering Rellich's question, E.~Wienholtz~\cite{Wienholtz} proved the following result:

\begin{thm}\label{T:Wien}
Let  $H_{g,b,V}$ be elliptic, $g^{jk}=g^{kj}$ are bounded and in $C^{3}(\RR^n)$,
$b_j\in C^3(\RR^n)$, $V$ real continuous function.
Suppose that $H_{g,b,V}$
is semibounded below on $C_c^{\infty}(\RR^n)$.
Then it is essentially self-adjoint.
\end{thm}
\noindent (Clearly neither Rellich nor Wienholtz were aware about the existence of
Povzner's paper. Wienholtz used a simpler method than Povzner. A simplified version of
Wienholtz result is explained in \cite{Glazman}.) Wienholtz~\cite{Wienholtz} also proved
the same statement in case when, instead of continuity, the potential $V$ is required to
belong to a global Stummel-type class.

The case of semi-boundedness of the operator $H_V$ is indeed essentially different from
the case of semi-bounded potential $V$. This is true even for $n=1$ as becomes clear from
an example by J.~Moser, which is described by F.~Rellich \cite{Rellich51} and also quoted
by H.~Kalf \cite{Kalf74}. This example provides a semi-bounded Schr\"odinger operator
$H_V$ in $L^2(\RR)$ with smooth $V$ such that $\Dom(H_{V,\max})$ is not contained in
$W^{1,2}(\RR)$ i.e. there exists $u\in \Dom(H_{V,\max})$ such that $u'\notin L^2(\RR)$.
This can not happen if $V$ is semi-bounded below as can be seen e.g. from Proposition
\ref{P:posquadform}.

There are many papers where the smoothness requirements on the coefficients in
\reft{Wien} were relaxed in different directions (see e.g. results quoted in
\cite{Kalf-75}). In particular, an important step was done by C.~Simader \cite{Simader78}
who considered semibounded below operator $H_V=-\Delta+V$  with $V=V_1+V_2$, where $0\leq
V_1\in\lloc^{2}(\RR^n)$ and $V_2$ satisfies a local Stummel type condition or a
local domination
condition. The proof is based on an observation that in this case
$\Dom(H_{V,\max})\subset W^{1,2}_{\loc}$ (which is the most difficult part of the proof),
and this is sufficient for the essential self-adjointness. {}Following \cite{sh}, we use
the Simader ideas in the geometric context of the present paper.

H.~Br\'ezis \cite{Brezis} introduced yet another
domination type local requirement on $V$ (``localized" self-adjointness)
which also implies the essential self-adjointness of
semi-bounded operator $H_V$ in $\RR^n$.
Note also that F.~S.~Rofe-Beketov and H.~Kalf~\cite{Kalf-Rofe-Beketov97}
unified Simader \cite{Simader78}  and Br\'ezis \cite{Brezis} results
by a refined use of localized self-adjointness.

The smoothness requirements for the metric coefficients $g^{jk}$  can be relaxed too.
{}For example, M.~A.~Perel'muter and Yu.~A.~Semenov \cite{Perel'muter-Semyonov} proved
essential self-adjointness  of $H_{g,0,V}$ under the conditions $V\ge 0$ and
$g^{jk}-\delta_{jk}\in W^{1,4}(\RR^n)$.

A Povzner-Wienholtz-type result for matrix-valued Sturm-Liouville operators
on intervals of $\R$, with coefficients in $L^1_{\loc}$, was established by
S.~Clark and F.~Gesztesy in their recent preprint \cite{Clark-Gesztesy}
which appeared after our paper was submitted.

Let us comment on the completeness condition in \reft{bound}. This condition can not be omitted.
In fact, N.~N.~Ural'ceva~\cite{Uralceva} and  S.~A.~Laptev~\cite{Laptev}
showed that there exist elliptic operators $H_{g,0,0}$ on $L^{2}(\RR^n)$, $n\ge 3$,
i.e. operators of the form
$$
\sum_{j,k=1}^n\frac{\pa}{\pa x^j}\left(g^{jk}(x)\frac{\pa}{\pa x^k}\right),
$$
with smooth positive definite $g^{jk}$ such that $C_c^\infty(\RR^n)$ is not dense
in the maximal domain $Q_{\max}$ of the corresponding quadratic form. If this is true then
the operator $H_{g,0,0}$ is not essentially self-adjoint.
This happens due to ``rapid growth" of $g^{jk}$,
so   the elements of the inverse matrix $g_{jk}$ are ``rapidly decaying" which means
that $\RR^n$ with the metric $g_{jk}$ is not complete.
(S.~A.~Laptev~\cite{Laptev} also gives some sufficient conditions
which insure that $C_c^\infty(\RR^n)$ is  dense
in  $Q_{\max}$.)

Amazingly,  V.~G.~Maz'ya \cite{Mazya67} (see also\cite[Sect.2.7]{Mazya-book})
established that the cases
$n=1,2$ are special here.   Using capacity, he proved that  $C_c^\infty(\RR^n)$ is always
dense in $Q_{\max}$ in case if $n=1$  or $n=2$.

It is possible that incompleteness of the metric  is compensated by a specific behavior
of $V$, so that the operator $H_{g,b,V}$ is self-adjoint even though the metric $g$ is
not complete (see e.g. A.~G.~Brusentsev~\cite{Brusentsev, Brusentsev-98},
H.~O.~Cordes~\cite[Ch.4]{Cordes2} and references there).

Let us make some comments on Schr\"odinger-type operators which may be not semibounded
below.  The first result as in \reft{main} (but on $M=\RR^n$ with standard metric and
measure, and with $q=q(|x|)$) for $H_V=-\Delta+V$, where $V\in\lloc^{\infty}(\RR^n)$, is
due to D.~B.~Sears~\cite{Sears} (see also \cite[Sect. 22.15 and 22.16]{Titchmarsh},
\cite[Ch. 3, Theorem 1.1]{bs}), who followed  an idea of an earlier paper by
E.~C.~Titchmarsh~\cite{Titchmarsh49}. {}F.~S.~Rofe-Beketov~\cite{RB} was apparently the
first to allow the minorant $q(x)$ to be not radially symmetric, though he did not
formulate the conditions geometrically.  He also proved, that the local inequality $V\geq
-q$ can be replaced by an operator inequality
\[
    H_V\ \geq \ -\delta\Delta-q(x)
\]
with a constant $\delta>0$.  This allows some potentials which are not bounded below even
locally. I.~Oleinik~\cite{ol, Oleinik94, Oleinik99} demonstrated  that this extends to
the case of manifolds. He also provided a geometric self-adjointness condition relating
it to the classical completeness in the situation which is not radially symmetric.
Oleinik's proof was  simplified by M.~Shubin \cite{sh3}, and then the result was extended
to magnetic Schr\"odinger operators in \cite{sh1}.

A first Sears-type result for $H_{b,V}$ with locally singular $V$ was obtained
by T.~Ikebe and T.~Kato~\cite{ik} (still with a radial minorant for $V$).
The paper \cite{sh1}  extends this result to complete Riemannian  manifolds, allows  non-radial
minorant of $V$ with Oleinik-type completeness condition, though requires $V$
to be locally bounded. By now we can also allow  appropriate
singularities of  $V$ for magnetic Schr\"odinger operators considered in \cite{sh1}
since such operators are a particular case of the ones allowed in Theorem \ref{T:main}.

A remarkable Sears-type result by A.~Iwatsuka \cite{Iwatsuka} seems to be the
only one which explicitly takes into account a direct involvement of magnetic field
in essential self-adjointness. In his result increasing magnetic field at infinity allows faster
fall off of the scalar potential $V$ to $-\infty$. (This fall off can be in fact arbitrarily fast
depending on the growth of the magnetic field.)

B.~M.~Levitan \cite{Levitan} suggested a new proof of the Sears theorem, based on consideration
of the hyperbolic Cauchy problem
\begin{equation}\label{E:cauchy}
\frac{\pa^2u}{\pa t^2}+Hu=0, \qquad u|_{t=0}=f, \quad \frac{\pa u}{\pa t}\Big|_{t=0}=g,
\end{equation}
where $H$ is the symmetric differential operator which we are investigating
(e.g. $H=H_{g,b,V}$ as in \eqref{E:oper}). Let us assume for a moment that $H$ is non-negative
(or semi-bounded below) on $C_c^\infty(\RR^n)$. Then essential self-adjointness of $H$
would follow if we know that for any $f,g\in C_c^\infty(\RR^n)$ the problem \eqref{E:cauchy}
has a unique solution which
is in $L^2(\RR^n_x)$ for every $t\in\RR$ (or even $t\in [0,t_0]$ with some $t_0>0$),
because different self-adjoint  extensions would generally produce different solutions
by the spectral theorem. The uniqueness in turn follows (by the Holmgren principle) from
the existence, provided a  globally finite propagation speed is secured, i.e. if
for every $t_0>0$, every bounded open set $G\subset\RR^n$ and every
$f,g\in C_{c}^{\infty}(G)$  there exists a bounded open set $\Omega\supset G$
such that the problem~\refe{cauchy} has a
solution $u(t,x)\in C^{\infty}([0,t_0]\times \RR^n)$ such that for all
$t\in[0,t_0]$, $\supp u(t,x)\subset\Omega$. A modification of this argument allows to consider
non-semi-bounded  operators as well. Establishing globally finite propagation speed of \refe{cauchy}
is a key ingredient in proving the essential self-adjointness.

A good explanation of the  abstract hyperbolic equation approach can be found in
Yu.~M.~Be\-re\-zansky \cite[Ch.VI, Sect.1.7]{Berezanskii-book}.
More details about the hyperbolic equation method can be found in the paper by
Yu.~B.~Orochko \cite{Orochko88} which also contains a good review with many
relevant references.

Yu.~M.~Berezansky's book \cite{Berezanskii-book} also contains an extensive discussion
on self-adjointness of operators generated by boundary value problems for elliptic and more
general operators.

Up to now we assumed at least that $V\in\lloc^{2}(\RR^n)$ which makes the minimal operator
defined on $C_c^\infty(\RR^n)$. However, T.~Kato~\cite{Kato74} pointed out that
if $V\in\lloc^{1}(\RR^n)$, one can still consider the maximal operator $H_{V,\max}$ associated
to $H_V=-\Delta+V$ as an operator
with the domain
$$
\Dom(H_{V,\max})=\{u\in L^{2}(\RR^n): Vu\in\lloc^{1}(\RR^n)
\text{ and }H_Vu\in L^{2}(\RR^n)\},
$$
where $H_Vu$ is a priori defined as a distribution.
The question which can be asked then is whether  the operator $H_{V,\max}$ is self-adjoint.
Besides, minimal operator $H_{V,\min}$ is then defined as the restriction of
$H_{V,\max}$ to compactly supported elements of $\Dom(H_{V,\max})$.
Then we can ask whether it is true that $H_{V,\min}$ is densely defined
and  $H_{V,\min}^*=H_{V,\max}$. More general operators of the form \eqref{E:oper}
can also be considered. We refer the reader to I.~Knowles \cite{Knowles78},
M.~Faierman and I.~Knowles \cite{Faierman-Knowles},
H.~Kalf and F.~S.~Rofe-Beketov \cite{Kalf-Rofe-Beketov} and Yu.~B.~Orochko \cite{Orochko88}
for results and references  in this direction.
In particular, Orochko treated locally integrable potentials
by the hyperbolic equation method.

Let us discuss some essential self-adjointness results for operators on manifolds. The
essential self-adjointness of the Laplace-Beltrami operator on a complete Riemannian
manifold was established by M.~Gaffney \cite{ga} (not only in $L^2(M)$ but also in the
standard $L^2$-spaces of differential forms). H.~O.~Cordes \cite{Cordes} (see also
\cite[Ch.4]{Cordes2})  established essential self-adjointness of powers of
Schr\"odinger-type operators with positive potentials (in $L^2$ spaces defined by a
measure which is unrelated with the metric, as in the present paper). He did it not only
on complete Riemannian manifolds but also in case when non-completeness is compensated by
an appropriate behavior of the potential, e.g. in domains $G\subset M$ of a complete
Riemannian manifold $M$ if the potential $V(x)$ is estimated below by $\gamma(\dist(x,\pa
G))^{-2}$ with an appropriate constant $\gamma>0$. Cordes uses a ``stationary" approach.
P.~Chernoff \cite{Chernoff} used the hyperbolic equation approach to establish similar
results. (In a later paper \cite{Chernoff77} he extended his results to the case of
singular potentials $V$.) We refer the reader to H.~O.~Cordes' book \cite[Ch.4]{Cordes2}
for more results about these topics. T.~Kato \cite{Kato2} extended Chernoff's arguments
from \cite{Chernoff} to prove the essential self-adjointness of powers of $H_V$, where
$H_V\geq -a-b|x|^2$ and $V$ is smooth  (in case $H_V=-\Delta+V$ on $M=\RR^n$ with the
standard metric and measure). A.~A.~Chumak \cite{Chumak} used the hyperbolic equation
method to prove essential self-adjointness of semi-bounded below operators $H_{g,0,V}$ on
complete Riemannian manifolds.

Let us also mention a subtle result of H.~Donnelly and N.~Garofalo~\cite{Donnelly-Garofalo}
who proved essential self-adjointness for $H_V=-\Delta+V$ on a complete Riemannian manifold
under the requirement that $V$ is locally in a Stummel-type class
except at a point $x_0$, and has a negative minorant $Q(x)$ such that
$Q(x)=\beta_0\dist(x,x_0)^{-2}$ near $x_0$ (which exactly matches the borderline case
in  Example  \ref{SS:ex}, so it is practically optimal), and near infinity
$Q(x)=-c\dist(x,x_0)^2$ with a constant $c>0$.

{}For some other related results on Schr\"odinger-type operators on manifolds see e.g.
references given in \cite{Cordes2, sh3, sh1, sh}. For higher order operators in $\RR^n$
see M.~Schechter \cite{Schechter} and F~.S.~Rofe-Beketov \cite{Rofe-Beketov3}. For higher
order operators on manifolds of bounded geometry see \cite{Shubin92}.


\end{document}